\documentclass[a4paper,11pt]{amsart}
\usepackage{comment}

\usepackage{lipsum}
\usepackage{url}

\usepackage{hyperref}

\usepackage[english]{babel}
\usepackage[T1]{fontenc}

\usepackage{a4wide}

\usepackage[utf8]{inputenc}
\usepackage{amsmath,amssymb,mathrsfs}
\usepackage{amscd}
\usepackage{pgf,tikz}
\usetikzlibrary{arrows}
\usetikzlibrary{patterns}
\numberwithin{equation}{section}
\newtheorem{theo}{Theorem}[section]
\newtheorem{cor}{Corollary}[section]
\newtheorem{cnj}{Conjecture}[section]
\newtheorem{prop}{Proposition}[section]
\newtheorem{lem}{Lemma}[section]

\newtheorem{de}{Definition}[section]
\newtheorem{rem}{Remark}[section]
\numberwithin{figure}{section}

\newcommand{\chh}{{H}}
\newcommand{\chho}{{H_0}}
\newcommand{\bn}{{\mathbb{N}}}
\newcommand{\bz}{{\mathbb{Z}}}
\newcommand{\br}{{\mathbb{R}}}
\newcommand{\bc}{{\mathbb{C}}}
\newcommand{\spp}{{\mathcal{S}_p}}
\newcommand{\sd}{{\mathcal{S}_2}}
\newcommand{\sinf}{{\mathcal{S}_\infty}}
\newcommand{\one}{{\bf 1}}

\newcommand{\xp}{X_\perp}

\title{On eigenvalue accumulation for non-self-adjoint magnetic operators}

\author{Diomba \textsc{Sambou}}
\address{Departamento de Matem\'aticas, Facultad de Matem\'aticas, 
Pontificia Universidad Cat\'olica de Chile, Vicu\~na Mackenna 
4860, Santiago de Chile}

\email{disambou@mat.uc.cl}

\thanks{The author is partially supported by the Chilean 
Program \textit{N\'ucleo Milenio de F\'isica Matem\'atica
RC$120002$}. The author gratefully acknowledges the many 
helpful suggestions of V. Bruneau during the preparation of the 
paper.}


\begin{document}

\begin{abstract}
In this work, we use regularized determinants to 
study the discrete spectrum generated by relatively compact 
non-self-adjoint perturbations of the magnetic Schrödinger 
operator $(-i\nabla - \textbf{\textup{A}})^{2} - b$ in $\br^3$, 
with constant magnetic field of strength $b>0$. The distribution 
of the above discrete spectrum near the Landau levels $2bq$, 
$q \in  \mathbb{N}$, is more interesting since they
play the role of thresholds of the spectrum of the free operator. 
First, we obtain sharp upper bounds on the number of 
complex eigenvalues near the Landau levels. 
Under appropriate hypothesis, we then prove the presence of 
 an infinite number of complex eigenvalues near each Landau 
level $2bq$, $q \in \mathbb{N}$, and the existence of sectors 
free of complex eigenvalues. We  also prove that the eigenvalues 
are localized in certain sectors adjoining the Landau levels. 
In particular, we provide an adequate answer to the open problem 
from \cite{dio} about the existence of complex eigenvalues 
accumulating near the Landau levels. 
Furthermore, we prove that the Landau levels are the only 
possible accumulation points of the complex eigenvalues.
\end{abstract}

\maketitle

\section{Introduction and motivations}\label{s1}

Presently, there is an increasing interest of mathematical physics community in the 
spectral theory of non-self-adjoint differential operators. 
Several results on the
discrete spectrum generated by non-self-adjoint perturbations have been established for the quantum Hamiltonians. Still, most 
of them give Lieb-Thirring type inequalities or upper bounds on certain distributional characteristics of eigenvalues, 
\cite{fra,bru1,bor,dem,demu,han,gol,wan,cue,dio,dub}
\big(for an extensive reference list on the subject, see for instance the 
references given in \cite{wan,cue}\big). Otherwise, results on spectral properties 
on non-self-adjoint operators can be found in the article by Sjöstrand \cite{sjo} 
and the references given there. In most of the above papers, the non-trivial 
question of the existence of complex eigenvalues near the essential spectrum
is not treated and stays open. 

For instance, in \cite{wan}, Wang studied $-\Delta + V$ in 
$L^2 \big( \br^n \big)$, $n \geq 2$, where the potential $V$ is dissipative. 
That is, 
\begin{equation}
V(x) = V_{1}(x) - iV_{2}(x),
\end{equation}
where $V_{1}$ and $V_{2}$ are two measurable functions such that $V_{2}(x) 
\geq 0$, and $V_{2}(x) > 0$ on an open non empty set. He showed that if the 
potential decays faster than $|x|^{-2}$, then the origin is not an 
accumulation point of the complex eigenvalues. For more general complex 
potentials without sign restriction on the imaginary part, it is still 
unknown whether the origin can be an accumulation point of complex eigenvalues 
or not. In this connection, the author \cite{dio1} proves  the existence of 
complex eigenvalues near the Landau levels together with their localization property
for non-self-adjoint two-dimensional Schrödinger operators with constant 
magnetic field. 

Motivated by Wang's work \cite{wan}, the current paper is devoted to 
the same type of results on eigenvalues near the Landau levels for the 
three-dimensional Schrödinger operator with constant magnetic 
field. Now, the essential spectrum of the operator under consideration 
equals $\br_+$, and the Landau levels play the role of thresholds. 
Consequently, the situation is more complicated than 
the non-self-adjoint case of the two-dimensional Schrödinger operator studied in 
\cite{dio1}, where the essential spectrum coincides with the (discrete) 
set of the Landau levels.

The magnetic field $\textbf{B}$ is generated by the magnetic potential 
$\textbf{A} = (-\frac{bx_2}2, \frac{bx_1}2, 0)$. Namely, 
$\textbf{B} = \text{curl} \hspace{0.1cm} \textbf{A} = (0,0,b)$
with constant direction, where $b > 0$ is a constant giving the
strength of the magnetic field. Then, the magnetic Schrödinger operator 
is defined by
\begin{equation}\label{op0}
H_0 := (-i\nabla - \textbf{A})^2  - b =  
\left( -i \frac{\partial}{\partial x_1} + \frac{b}{2} 
x_2 \right)^2 + \left( -i \frac{\partial}{\partial x_2} - 
\frac{b}{2} x_1 \right)^2  + \left( -i \frac{\partial}{\partial x_3} 
\right)^2 - b,
\end{equation}
in $L^2(\br^3)$ with $x = (x_1,x_2,x_3) \in \br^{3}$. Actually, $H_0$ 
is the self-adjoint operator associated with the closure $\overline{q}$ 
of the quadratic form 
\begin{equation}
q(u) = \int_{\br^{3}} \left( \big\vert (-i\nabla - \textbf{A}) u(x) \big\vert^{2} 
 - b \vert u(x) \vert^{2} \right) dx,
\end{equation}
originally defined on $C_{0}^{\infty}(\br^{3})$. The form domain 
$D(\overline{q})$ of $\overline{q}$ is the magnetic Sobolev space
$
H_{\textbf{A}}^{1}(\br^{3}) := \big \lbrace u \in L^{2}(\br^{3}) : 
(-i\nabla - \textbf{A}) u \in L^{2}(\br^{3}) \big\rbrace
$,
\big(see for instance \cite{lt2}\big). Setting $\xp := (x_1, x_2) 
\in \br^2$ and $L^2(\br^3) 
= L^2(\br_{\xp}^2) \otimes L^2(\br_{x_3})$, $H_0$ can be rewritten in the form
\begin{equation}\label{opll}
H_0 = H_{\rm Landau} \otimes I_3 + I_{\perp} \otimes
\left( - \frac{\partial^2}{\partial x_{3}^2} \right).
\end{equation}
Here,
\begin{equation}\label{opl}
H_{\text{Landau}}: = \left( -i \frac{\partial}{\partial x_1} + \frac{b}{2} 
x_2 \right)^2 + \left( -i \frac{\partial}{\partial x_2} - 
\frac{b}{2} x_1 \right)^2 - b
\end{equation}
is the shifted Landau Hamiltonian, self-adjoint in $L^2(\br^2)$, and 
$I_3$, $I_\perp$ are the identity operators in $L^2(\br_{x_3})$ and 
$L^2(\br_{\xp})$ respectively. It is well known \big(see for instance 
\cite{avr,dim}\big) that the spectrum of $H_{\text{Landau}}$ 
consists of the so-called Landau levels $\Lambda_q := 2bq$, 
$q \in \bn : = \{0,1,2, \ldots\}$, and 
${\rm dim}\,{\rm Ker}(H_{\text{Landau}} - \Lambda_q) = \infty$. Hence,
\begin{equation*}
\sigma ( H_0 ) = \sigma_{\rm ac} ( H_0 ) = [0, +\infty),
\end{equation*}
and, once again, the Landau levels play the role of thresholds of this spectrum. 

\begin{rem}
\textup{Looking at \eqref{opll} as well as the structure of the spectrum of 
$H_{\text{Landau}}$ and the one of $-\frac{\partial^2}{\partial x_3^2}$, one sees 
that the structure of $H_0$ is quite close to the one of the (free) quantum waveguide 
Hamiltonians.}
\end{rem}

Let us introduce some important definitions. Let $M$ be a closed linear operator 
acting on a separable Hilbert space $\mathscr{H}$. If $z$ is an isolated 
point of $\sigma(M)$, the spectrum of $M$, let $\gamma$ be a small positively 
oriented circle centred at $z$ and containing $z$ as the only point of 
$\sigma(M)$. 

\begin{de}[Discrete eigenvalue]
The point $z$ is said to be a discrete eigenvalue of
$M$ if its algebraic multiplicity is finite and
\begin{equation}\label{eq1,0}
\text{mult}(z) := \text{rank} \left( \frac{1}{2i\pi} 
\int_{\gamma} (M - \zeta)^{-1} d\zeta \right).
\end{equation}
\end{de}
Note that we have
$\text{mult}(z) \geq \text{dim} \big( \text{Ker} (M - z) \big)$, the geometric 
multiplicity of $z$. The inequality becomes an equality if $M$ is self-adjoint. 

\begin{de}
[Discrete spectrum]
The discrete spectrum of $M$ is defined by 
\begin{equation}\label{eq1,1}
\sigma_{\textup{\textbf{disc}}}(M) := \big\lbrace z \in \bc : z 
\hspace*{0.1cm} \textup{\textit{is a discrete eigenvalue of} $M$} 
\big\rbrace.
\end{equation}
\end{de}

\begin{de}
[Essential spectrum]
The essential spectrum of $M$ is defined by
\begin{equation}\label{eq1,2}
\sigma_{\textup{\textbf{ess}}}(M) := \big\lbrace z \in \bc : 
\textup{$M - z$ \textit{is not a Fredholm operator}} \big\rbrace.
\end{equation} 
It is a closed subset of $\sigma(M)$.
\end{de}

The purpose of this paper is to investigate the distribution of the 
discrete spectrum near the essential spectrum of the perturbed 
operator 
\begin{equation}\label{eq1,3}
H := \chho + W \qquad \text{on} \qquad Dom(\chho),
\end{equation}
where $W : \br^3 \longrightarrow \bc$ is a non-self-adjoint relatively compact 
perturbation with respect to $H_0$. In \eqref{eq1,3}, $W$ is identified with 
the multiplication operator by the function (also denoted) $W$. In the sequel, $W$ is supposed 
to satisfy some general assumptions  \big(see \eqref{eq1,7} \big).

To put our results in perspective, let us first discuss known results in the case of self-adjoint perturbations. 
It is well known \big(see for instance \cite[Theorem 1.5]{avr}\big) that 
if $W : \br^{3} \longrightarrow \br$ satisfies
\begin{equation}\label{eq1,4}
W(x) \leq -C \one_U(x), \qquad x \in \br^3,
\end{equation}
for some constant $C > 0$ and some non-empty open set $U \subset \br^3$, then the 
discrete spectrum of $H$ is infinite. Moreover, if $W$ is axisymmetric ($i.e.$ depends 
only on $\vert \xp \vert$ and $x_3$) and verifies \eqref{eq1,4}, then it is known 
\big(see for instance \cite[Theorem 1.5]{avr}\big) that $H$ has an infinite number 
of eigenvalues embedded in the essential spectrum. In the case where $W$ is axisymmetric verifying
\begin{equation}\label{eq1,5}
W(x) \leq -C \one_S(\xp) (1 + \vert x_3 \vert)^{-m_3},
\qquad m_3 \in (0,2), \qquad x = (\xp,x_3) \in \br^3,
\end{equation}
for some constant $C > 0$ and some non-empty open set $S \subset \br^2$, it is also proved 
\big(see \cite{rag,rage}\big) that below each Landau level $2bq, q \in \bn$, there is
an infinite sequence of discrete eigenvalues of $H$ converging to $2bq$. In \cite{bon, 
bo}, the resonances of the operator $H$ near the Landau levels have been investigated for 
self-adjoint potentials $W$ decaying exponentially in the direction of the magnetic 
field. Namely,
\begin{equation}\label{eq1,51}
W(x) = \mathcal{O} \big( (1 + \vert \xp \vert)^{-m_\perp}
\exp \big( -N \vert x_3 \vert \big), \qquad x = (\xp,x_3) \in \br^3,
\end{equation} 
with $m_\perp > 0$ and $N > 0$. 
Other results on the distribution  of discrete spectrum for magnetic quantum Hamiltonians 
perturbed by self-adjoint electric potentials can be found in \cite[Chap. 11-12]{iv}, \cite{pus,ra,rai,mel,sob,tam,roz} and the references therein.

Throughout this paper, our minimal assumption on $W$ defined by
\eqref{eq1,3} is the following:
\begin{equation}\label{eq1,7}
\textit{\textbf{Assumption (A1):}} 
\begin{cases}
\bullet \hspace{0.6mm} W \in L^\infty (\br^3,\bc), W(x) 
= \mathcal{O} \big( F(\xp) G(x_3) \big), x = (\xp,x_3) \in \br^3, \\
\bullet \hspace{0.6mm} F \in \bigl( L^\frac{p}{2} \cap L^\infty 
\bigr) \big( \br^{2},\br_{+}^{\ast} \big) \hspace{0.4mm} 
\textup{for some $p \geq 2$}, \\
\bullet \hspace{0.6mm} \br_{+}^{\ast} \ni G(x_3) = \mathcal{O}
\big( \langle x_3 \rangle^{-m} \big), m > 3,
\end{cases}
\end{equation}
where $\langle y \rangle := \sqrt{1 + \vert y \vert^2}$ for
$y \in \br^d$.

\begin{rem}
\textup{Typical example of potentials satisfying \textit{Assumption (A1)} is the 
special case of potentials $W : \br^3 \rightarrow \bc$ such that
\begin{equation}\label{eq1,81}
W(x) = \mathcal{O} \big( \langle \xp \rangle^{-m_{\perp}} 
\langle x_3 \rangle^{-m} \big), \quad m_{\perp} > 0, \quad  m > 3.
\end{equation}
We can also consider the class of potentials $W : \br^3 \rightarrow \bc$ such that
\begin{equation}\label{eq1,82}
W(x) = \mathcal{O} \big( \langle \textup{\textbf{x}} \rangle^{-\alpha} 
\big), \quad \alpha > 3.
\end{equation}
Indeed, condition \eqref{eq1,82} implies that \eqref{eq1,81} holds for any $m \in 
(3,\alpha)$ and $m_{\perp} = \alpha - m > 0$.}
\end{rem}

Under \textit{Assumption (A1)}, we establish (see Lemma \ref{lem3,1}) that the weighted 
resolvent $\vert W \vert^{\frac{1}{2}} (\chho - z)^{-1}$ belongs to the Schatten-von 
Neumann class $\spp$ (see Subsection \ref{ss3.1} where the classes $\spp$, $p \geq 1$ 
are introduced). Consequently, $W$ is relatively compact with respect to $\chho$. Then, 
from Weyl's criterion on the invariance of the essential spectrum, it follows 
that
\begin{equation}
\sigma_{\text{\textbf{ess}}}(\chh) = 
\sigma_{\text{\textbf{ess}}}(\chho) = [0,+\infty).
\end{equation}
However, the electric potential $W$ may generate (complex) discrete eigenvalues \big($\sigma_{\text{\textbf{disc}}}(\chh)$\big) that can only accumulate to 
$\sigma_{\text{\textbf{ess}}}(\chh)$, see \cite[Theorem 2.1, p. 373]{goh}. A natural 
question is to sharpen the rate of this accumulation by studying the distribution of $\sigma_{\text{\textbf{disc}}}(\chh)$ near $[0,+\infty)$, in particular near the spectral 
thresholds $2bq, q \in \bn$. 
Motivated by this problem,  the following 
result  \cite{dio}, often called a generalized Lieb-Thirring type inequality \big(see Lieb-Thirring 
\cite{lt1} for original work\big), is obtained by using complex analysis tools developed 
by Borichev-Golinskii-Kupin \cite{bor}.

\begin{theo}{\cite[Theorem 1.1]{dio}}\label{theo1}

\noindent
Let $\chh := \chho + W$ with $W : \br^3 \longrightarrow \bc$ being bounded 
and satisfying the inequality $W(x) = \mathcal{O} \big( F(x) G(x_3) \big)$, 
where $F \in \big( L^\infty \cap L^{p} \big) (\br^3)$, $p \geq 2$, and 
$G \in \big( L^\infty \cap L^2 \big)(\br)$. Then, for any 
$0 < \varepsilon < 1$, we have
\begin{equation}\label{eq1,9}
\small{\displaystyle \sum_{z \hspace{0.5mm} \in \hspace{0.5mm} 
\sigma_{\textup{\textbf{disc}}}(H)} 
\frac{\textup{dist} \big( z,[\Lambda_{0},+\infty) \big)^{\frac{p}{2} 
+ 1 + \varepsilon} 
\hspace{0.5mm} \textup{dist} \big( z,\cup_{q=0}^{\infty} 
\lbrace \Lambda _{q} \rbrace \big)^{(\frac{p}{4} - 1 + \varepsilon)_{+}}}{(1 + 
\vert z \vert)^{\gamma}} \leq C_{1} K},
\end{equation}
where $\gamma > d + \frac{3}{2}$, $C_{1} = C(p,b,d,\varepsilon)$ and
$$
K := \Vert F \Vert_{L^{p}}^{p} \big{(} \Vert G \Vert_{L^{2}} + 
\Vert G \Vert_{L^{\infty}} \big{)}^{p} 
\big{(} 1 + \Vert W \Vert_{\infty} \big{)}^{d  + \frac{p}{2} + 
\frac{3}{2} + \varepsilon}.
$$ 
Here, $r_{+} := \max \hspace{0.5mm} (r,0)$ for $r \in \br$.
\end{theo}

We give few comments on the above result to make the connection with the present problem more explicit. Let
$(z_{\ell})_{\ell} \subset \sigma_{\textup{\textbf{disc}}}(\chh)$ 
be a sequence of complex eigenvalues that converges  non-tangentially to a Landau level $\Lambda_q = 
2bq, q \in \bn$. Namely,
\begin{equation}
\vert \Re(z_\ell) - 2bq \vert \leq C \hspace{0.5mm} 
\vert \Im(z_\ell) \vert,
\end{equation} 
for some constant $C > 0$ \big(\textbf{(iii)} of Theorem \ref{theo4} implies that a 
such sequence exits if $W$ in Theorem \ref{theo1} satisfies the required conditions\big). 
Thus, bound \eqref{eq1,9} implies (taking  a subsequence if necessary), that
\begin{equation}\label{eq1,10}
\sum_{\ell} \textup{dist} \big( z_{\ell},\cup_{q=0}^{\infty} 
\lbrace \Lambda_{q} \rbrace \big)^{(\frac{p}{2} + 1 + \varepsilon) + 
(\frac{p}{4} - 1 + \varepsilon)_{+}} 
< \infty.
\end{equation}
Formally, \eqref{eq1,10} means that the sequence $(z_\ell)_\ell$ converges to the Landau 
level with a rate convergence larger than $\frac{1}{(\frac{p}{2} + 1 + \varepsilon) + 
(\frac{p}{4} - 1 + \varepsilon)_{+}}$. This means that the ``convergence exponent'' of such sequences 
near the Landau levels is a monotone function of $p$. However, even if 
Theorem \ref{theo1} allows to estimate formally the rate accumulation of the complex 
eigenvalues (near the Landau levels), it does not prove their existence.

Two important points of the present paper are to be taken into account. First, we prove the presence
of infinite number of complex eigenvalues of $\chh$ near each Landau level $2bq$, $q \in 
\bn$, for certain classes of potentials $W$ satisfying \textit{Assumption (A1)}. Second, 
we prove that the Landau levels are the only possible accumulation points of the discrete eigenvalues, see Theorem \ref{theo6}. It is worth mentioning that we expect this to be a general phenomenon. 

Our techniques are close to those from \cite{bon} used for the study of the resonances 
near the Landau levels for self-adjoint electric potentials. Firstly, we obtain sharp 
upper bound on the number of discrete eigenvalues in small annulus around a Landau 
level $2bq$, $q \in \bn$, for general complex potentials $W$ satisfying \textit{Assumption 
(A1)} (see Theorem \ref{theo2}). Secondly, under appropriate assumption \big(see 
\textit{Assumption (A2)} given by \eqref{eq2,5}\big), we obtain a special upper 
bound on the number of discrete eigenvalues outside a semi-axis in annulus centred 
at a Landau level (see Theorem \ref{theo3}). Under additional hypothesis, \big(see \textit{Assumption (A3)} given by \eqref{eq2,7}\big), we establish corresponding lower 
bounds implying the existence of an infinite number of discrete eigenvalues or the absence 
of discrete eigenvalues in some sectors adjoining the Landau levels $2bq$, $q \in 
\bn$ (see Theorem \ref{theo4}). In particular, we derive from Theorem \ref{theo4} a 
criterion of non-accumulation of complex eigenvalues of $H$ near the Landau levels, 
see Corollary \ref{c1} (see also Conjecture \ref{cj2,1}). Loosely speaking, our methods 
can be viewed as a Birman-Schwinger principle applied to the non-self-adjoint perturbed 
operator $\chh$ (see Proposition \ref{prop3,2}). By this way, we reduce the study of the 
discrete eigenvalues near the essential spectrum, to the analysis of the zeros of a holomorphic 
regularized determinant.

The paper is organized as follows. Section \ref{s2} is devoted to the statement of our 
main results. In Section \ref{s3}, we recall useful properties on regularized determinant 
defined for operators lying in the Schatten-von Neumann classes $\spp$, $p \geq 1$. Furthermore,  we establish a first reduction of the study of the complex eigenvalues in a neighbourhood 
of a fixed Landau level $2bq$, $q \in \bn$, to that of the zeros of a holomorphic function. 
In Section \ref{s4}, we establish a decomposition of the weighted free resolvent, which is 
crucial for the proofs of our main results in Sections \ref{s5}-\ref{s7}. Section \ref{sa} 
is a brief Appendix presenting tools on the index of a finite meromorphic operator-valued 
function.

\section{Formulation of the main results}\label{s2}

We start this section with a list of useful notations and definitions.% are needed. 

We denote $P_q$ the orthogonal projection onto ${\rm Ker}(H_{\text{Landau}} - 
\Lambda_q)$, $\Lambda_q = 2bq$, $q \in \bn$. 

For $W$ satisfying \textit{Assumption (A1)}, introduce $\textbf{\textup{W}}$ the 
multiplication operator by the function (also denoted) $\textbf{\textup{W}} : \br^2 \longrightarrow 
\br$ defined by 
\begin{equation}\label{eq2,1}
\displaystyle \textbf{\textup{W}}(\xp) := \frac{1}{2}
\int_\br \vert W(\xp,x_3) \vert dx_3.
\end{equation}
By \cite[Lemma 5.1]{ra}, if $U \in L^p(\br^2)$, $p \geq 1$, 
then $P_qUP_q \in \spp$ for any $q \in \bn$. According to \eqref{eq1,7},
$\textbf{\textup{W}}(\xp) = \mathcal{O} \big( F(\xp) \big) =
\mathcal{O} \big( F^{\frac{1}{2}}(\xp) \big) $. Thus, since $F^{\frac{1}{2}} \in L^p 
(\br^2)$, the Toeplitz operator $P_q \textbf{\textup{W}} P_q \in \spp$ for any 
$q \in \bn$. 

Our results are closely related to the quantity $\textup{Tr} \hspace{0.4mm} 
\one_{(r,\infty)} \big( P_q \textbf{\textup{W}} P_q \big)$, $r > 0$. When the function $\textbf{\textup{W}} = U$ admits a power-like decay, a exponential decay, or is compactly 
supported, then asymptotic expansions of 
$\textup{Tr} \hspace{0.4mm} \one_{(r,\infty)} \big( P_q \textbf{\textup{W}} P_q \big)$ 
as $r \searrow 0$ are well known:

\textbf{(i)} If $0 \leq U \in C^{1} (\mathbb{R}^{2})$ satisfies 
$U(\xp) = u_{0}\big(\xp / \vert \xp \vert\big) \vert \xp \vert^{-m} ( 1 + o(1) \big)$, 
$\vert \xp \vert \rightarrow \infty$, $u_{0}$ being a non-negative continuous function 
on $\mathbb{S}^{1}$ not vanishing identically, and $\vert \nabla U(\xp) \vert \leq C_{1} 
\langle \xp \rangle^{-m-1}$ with some constants $m > 0$ and $C_{1} > 
0$, then by \cite[Theorem 2.6]{ra}
\begin{equation}\label{éq2,7}
\textup{Tr} \hspace{0.4mm} \one_{(r,\infty)} 
\big( P_{q} U P_{q} \big) = 
C_{m} r^{-2/m} \big( 1 + o(1) \big), \hspace{0.2cm} r \searrow 0,
\end{equation}
where $C_{m} := \frac{b}{4\pi} \int_{\mathbb{S}^{1}} u_{0}(t)^{2/m} dt$. Note that in 
\cite[Theorem 2.6]{ra}, \eqref{éq2,7} is stated in a more general version including higher 
even dimensions $n = 2d$, $d \geq 1$.

\textbf{(ii)} If $0 \leq U \in L^{\infty} (\mathbb{R}^{2})$ satisfies $\ln U(\xp) 
= -\mu \vert \xp \vert^{2\beta} \big( 1 + o(1) \big)$, $\vert x \vert \rightarrow \infty$, 
with some constants $\beta > 0$ and $\mu > 0$, then by \cite[Lemma 3.4]{raik}
\begin{equation}\label{éq2,8}
\textup{Tr} \hspace{0.4mm} \one_{(r,\infty)} 
\big( P_{q} U P_{q} \big) = 
\varphi_{\beta}(r) \big( 1 + o(1) \big), \hspace{0.2cm} r \searrow 0,
\end{equation}
where we set for $0 < r < \textup{e}^{-1}$
$$
\varphi_{\beta}(r) :=
 \begin{cases}
 \frac{1}{2} b \mu^{-1/\beta} \vert \ln r \vert^{1/\beta} & \text{if } 
 0 < \beta < 1,\\
 \frac{1}{\ln(1 + 2\mu/b)} \vert \ln r \vert & \text{if } \beta = 1,\\
 \frac{\beta}{\beta - 1} \big{(} \ln \vert \ln r \vert \big{)}^{-1} \vert 
 \ln r \vert & \text{if } \beta > 1.
 \end{cases}
$$

\textbf{(iii)} If $0 \leq U \in L^{\infty} (\mathbb{R}^{2})$ is compactly supported 
and if there exists a constant $C > 0$ such that $C \leq U$ on an open non-empty subset 
of $\mathbb{R}^{2}$, then by \cite[Lemma 3.5]{raik}
\begin{equation}\label{éq2,9}
\textup{Tr} \hspace{0.4mm} \one_{(r,\infty)} 
\big( P_{q} U P_{q} \big) = 
\varphi_{\infty}(r) \big( 1 + o(1) \big), 
\hspace{0.2cm} r \searrow 0,
\end{equation}
with
$\varphi_{\infty}(r) := \big( \ln \vert \ln r \vert 
\big{)}^{-1} \vert \ln r \vert, \hspace{0.2cm} 0 < r 
< \textup{e}^{-1}$.
Note that extensions of \cite[Lemmas 3.4 and 3.5]{raik} in higher even dimensions 
are established in \cite{mel}.

Now, introduce respectively the upper and lower half-planes by
\begin{equation}\label{eq2,2}
\bc_\pm := \big\lbrace z \in \bc : \pm \Im(z) > 0 \big\rbrace.
\end{equation}
For a fixed Landau level $\Lambda_{q} = 2bq$, $q \in \bn$, and $0 \leq a_1 < a_2 \leq 2b$, 
define the ring
\begin{equation}\label{eq2,3bb}
\Omega_{q}(a_1,a_2) := \big\lbrace z \in \bc: a_1 < \vert \Lambda_{q} - 
z \vert < a_2  \big\rbrace,
\end{equation}
and the half-rings
\begin{equation}\label{eq2,3b}
\Omega_{q}^{\pm}(a_1,a_2) := \Omega_{q}(a_1,a_2) \cap \bc_\pm.
\end{equation}
For $\nu > 0$, we introduce the domains
\begin{equation}\label{eq2,3bbb}
\Omega_{q,\nu}^{\pm}(a_1,a_2) := \Omega_{q}^{\pm}(a_1,a_2) \cap 
\big\lbrace z \in \bc : \vert \Im(z) \vert > \nu \big\rbrace.
\end{equation}
Our first main result gives an upper bound on the number of discrete eigenvalues 
in small half-rings around a Landau level $2bq$, $q \in \bn$.

\begin{theo}[Upper bound]\label{theo2}
Assume that \textit{Assumption (A1)} holds with $0 < \Vert W \Vert_\infty < 2b$ small enough. 
Then, there exists $0 < r_{0} < \sqrt{2b}$ such that for any $r > 0$ with
$r < r_{0} < \sqrt{\frac{5}{2}}r$ and any $q \in \bn$,
\begin{equation}\label{ub1}
\# \big\lbrace z \in \sigma_{\textup{\textbf{disc}}}(\chh) \cap \Omega_{q,\nu}^{\pm}(r^2,4r^2)
\big\rbrace = \mathcal{O} \Big( \textup{Tr} \hspace{0.4mm} \one_{(r,\infty)} \big( 
P_{q} \textbf{\textup{W}} P_{q} \big) \vert \ln r \vert \Big),
\end{equation}
$0 < \nu < 2r^2$. In particular, if the function $\textbf{\textup{W}}$ is compactly 
supported, then $\textup{Tr} \hspace{0.4mm} \one_{(r,\infty)} \big( P_{q} \textbf{\textup{W}} P_{q} \big) = \mathcal{O} \left( \big(\ln \vert \ln r \vert \big)^{-1} \vert \ln r \vert \right)$ 
as $r \searrow 0$.
\end{theo}

\noindent
In order to formulate the rest of our main results, it is necessary to make additional 
restrictions on $W$. Namely,
\begin{equation}\label{eq2,5}
\textit{\textbf{Assumption (A2):}} 
\begin{cases}
\textup{$W = e^{i\alpha} V$ with $\alpha \in \br \setminus \pi \bz$, and
$V : \text{Dom}(\chho) \longrightarrow L^2(\br^3)$ is the} \\
\textup{multiplication operator by the function 
$ V : \br^3 \longrightarrow \br$}.
\end{cases}
\end{equation}

\noindent
Note that in \textit{Assumption (A2)}, we can replace $e^{i\alpha}$ by any complex number 
$c = \vert c \vert e^{iArg(c)} \in \bc \setminus \br$. 

\begin{rem}\label{r2}
\textup{\textbf{(i)} In \eqref{eq2,5}, when $V$ is of definite sign ($i.e.$ $\pm V \geq 0$), 
since the change of the sign consists to replace $\alpha$ by $\alpha + \pi$, then it 
is enough to consider only $V \geq 0$.}

\textup{\textbf{(ii)} For $\pm \sin(\alpha) > 0$ and $V \geq 0$, the discrete eigenvalues 
$z$ of $H$ satisfy $\pm \Im(z) \geq 0$.} 
\end{rem}

\noindent
The next result gives an upper bound on the number of discrete eigenvalues outside 
a semi-axis, in small half-rings around a Landau level.

\begin{theo}[Upper bound, special case]\label{theo3}
Let $W$ satisfy Assumption (A1) with $0 < \Vert W \Vert_\infty < 2b$ small enough, 
$F \in L^1 (\br^2)$, and Assumption (A2) with $V \geq 0$, $\alpha = \pm \frac{3\pi}{4}$.
Then, for any $\theta > 0$ small enough, there exists $r_0 > 0$ such that for any 
$r > 0$ with $r < r_{0} < \sqrt{\frac{5}{2}}r$ and any $q \in \bn$,
\begin{equation}\label{ubs1}
\# \left\lbrace z \in \sigma_{\textup{\textbf{disc}}}(\chh) \cap 
\Omega_{q,\nu}^{\pm}(r^2,4r^2) \cap E_{q}^{\pm}(\alpha,\theta) \right\rbrace = 
\mathcal{O} \big( \vert \ln r \vert \big),
\end{equation}
$0 < \nu < 2r^2$, where $E_{q}^{\pm}(\alpha,\theta) := \Omega_{q}(0,2b) \setminus \left( 2bq + 
e^{i (2\alpha \mp \pi)} e^{i (-2\theta,2\theta)} (0,2b) \right)
$.
\end{theo}

\begin{rem}
\textup{Notice that in the setting $E_{q}^{\pm}(\alpha,\theta)$ above, we have just
excluded an angular sector of amplitude $4\theta$ around the semi-axis 
$z = 2bq + e^{i(2\alpha \mp \pi)} (0,2b)$.}
\end{rem}

\noindent
To get the existence of an infinite number of complex eigenvalues near the Landau 
levels, we need to assume at least that the function $\textbf{\textup{W}}$ defined by 
\eqref{eq2,1} has an exponential decay:
\begin{equation}\label{eq2,7}
\textit{\textbf{Assumption (A3):}}
\begin{cases}
\textbf{\textup{W}} \in L^{\infty}(\br^2), \quad 
\ln \textbf{\textup{W}}(\xp) \leq -C \langle \xp \rangle^2 \\
\textup{for some constant $C > 0$.} 
\end{cases}
\end{equation}

\begin{theo}[Sectors free of complex eigenvalues, upper 
and lower bounds]\label{theo4}
Under the assumptions and the notations of Theorem \ref{theo3} with the condition 
$F \in L^1 (\br^2)$ removed, for any $\theta > 0$ small enough, there exists 
$\varepsilon_0 > 0$ such that:

$\textup{\textbf{(i)}}$
For any $0 < \varepsilon \leq \varepsilon_0$, $\chh_\varepsilon := \chho + \varepsilon 
W$ has no discrete eigenvalues in
\begin{equation}\label{sec1}
\Omega_{q}^{\pm}(r^2,r_0^2) \cap E_{q}^{\pm}(\alpha,\theta), \quad r_0 \ll 1.
\end{equation}

$\textup{\textbf{(ii)}}$ If moreover $F \in L^1 (\br^2)$ in Assumption (A1), then 
there exists $r_0 > 0$ such that for any $0 < r <r_0$ and $0 < \varepsilon \leq 
\varepsilon_0$,
\begin{equation}\label{ubss1}
\small{ \# \left\lbrace z \in \sigma_{\textup{\textbf{disc}}}(\chh_\varepsilon) 
\cap \Omega_{q,\nu}^{\pm} \Big( \frac{4r^2}{9},\frac{9r^2}{4} \Big) \right\rbrace = 
\mathcal{O} \Big( \textup{Tr} \hspace{0.4mm} \one_
{(\frac{r}{2},4r)} \big( \varepsilon P_{q} \textbf{\textup{W}} P_{q} \big) \Big), \quad
0 < \nu < \frac{8r^2}{9}.}
\end{equation}

$\textup{\textbf{(iii)}}$ Let $\textbf{\textup{W}}$ satisfy Assumption (A3). Then, for 
any $0 < \varepsilon \leq \varepsilon_0$, there is an accumulation of discrete eigenvalues 
of $\chh_\varepsilon$ near $2bq$, $q \in \bn$, in a sector around the semi-axis 
$z = 2bq + e^{i(2\alpha \mp \pi)} ]0,+\infty)$, for
\begin{equation}\label{ca1}
\alpha \in \pm \left( \frac{\pi}{2},\pi \right).
\end{equation}
More precisely, there exists a decreasing sequence $(r_\ell)_\ell$ of positive numbers 
$r_\ell \searrow 0$ such that
\begin{equation}\label{lb1}
\small{ \# \left\lbrace z \in \sigma_{\textup{\textbf{disc}}}(\chh_\varepsilon) \cap 
\Omega_{q}^{\pm}(\varepsilon^2 r_{\ell + 1}^2,\varepsilon^2 r_{\ell}^2) \cap 
\left( 2bq + e^{i (2\alpha \mp \pi)} e^{i (-2\theta,2\theta)} (0,2b) \right) 
\right\rbrace \geq \textup{Tr} \hspace{0.4mm} \one_{(r_{\ell +1},r_\ell)} 
\big( P_{q} \textbf{\textup{W}} P_{q} \big).}
\end{equation}
\end{theo}

\begin{figure}[h]\label{fig 2}
\begin{center}
\tikzstyle{+grisEncadre}=[fill=gray!60]
\tikzstyle{blancEncadre}=[fill=white!100]
\tikzstyle{grisEncadre}=[fill=gray!25]
\tikzstyle{dEncadre}=[dotted]

\begin{tikzpicture}[scale=1]

\draw (0,0) -- (0:2) arc (0:180:2) -- cycle;

\node at (-0.25,-0.15) {\tiny{$2bq$}};
\node at (0.5,-0.16) {\tiny{$r^2$}};
\node at (2,-0.17) {\tiny{$r_0^2$}};

\draw [grisEncadre] (0,0) -- (0:2) arc (0:180:2) -- cycle;

%\draw [blancEncadre] (0,0) -- (0:0.5) arc (0:90:0.5) -- cycle;

\draw [+grisEncadre] (0,0) -- (39:2) arc (39:73:2) -- cycle;

\draw (0,0) -- (0:0.5) arc (0:180:0.5) -- cycle;

\draw (0,0) -- (-2.5,2.5) -- cycle;
\draw (-2.6,2.35) node[above] {$e^{i\alpha}\br_+$};
\draw (-2.45,0.6) node[above] {$\pi - \alpha$};
\draw (0,0) -- (135:2.05) arc (135:180:2.05) -- cycle;

\draw[->] [thick] (-3,0) -- (3,0);
\draw (3,0) node[right] {\tiny{$\Re(z)$}};

\draw[->] [thick] (0,-0.5) -- (0,2.5);
\draw (0,2.5) node[above] {\tiny{$\Im(z)$}};

\draw (0,0) -- (1.7,2.5);
\draw (2.5,2.4) node[above] {\tiny{$y = \tan (2\alpha - \pi) 
\hspace{0.5mm} (x - 2bq)$}};

\draw (0.95,1.75) node[above] {\tiny{$2\theta$}};
\draw (1.45,1.42) node[above] {\tiny{$2\theta$}};

\draw [dEncadre] [->] (-4,0.8) -- (-0.8,1.3);
\draw [dEncadre] [->] (-4,0.8) -- (1.1,0.35);
\draw (-4,0.8) node[left] {\tiny{$\Omega_q^+(r^2,r_0^2) \cap E_q^+(\alpha,\theta)$}};

\node at (0.45,0.5) {\tiny{$\times$}};
\node at (0.45,0.65) {\tiny{$\times$}};
\node at (0.63,0.9) {\tiny{$\times$}};
\node at (0.76,1.1) {\tiny{$\times$}};
\node at (0.8,0.95) {\tiny{$\times$}};
\node at (0.62,0.75) {\tiny{$\times$}};

\node at (0.3,0.6) {\tiny{$\times$}};
\node at (0.4,0.8) {\tiny{$\times$}};
\node at (0.6,1.1) {\tiny{$\times$}};
\node at (0.5,0.95) {\tiny{$\times$}};

\node at (0.35,0.35) {\tiny{$\times$}};
\node at (0.2,0.3) {\tiny{$\times$}};
\node at (0.25,0.45) {\tiny{$\times$}};
\node at (0.6,0.6) {\tiny{$\times$}};
\node at (0.8,0.8) {\tiny{$\times$}};
\node at (1.2,1.4) {\tiny{$\times$}};
\node at (1,1.6) {\tiny{$\times$}};
\node at (0.7,1.7) {\tiny{$\times$}};
\node at (1,1.2) {\tiny{$\times$}};
\node at (1.1,1) {\tiny{$\times$}};

\node at (0,4) {$W = e^{i\alpha} V$};
\node at (0,4.48) {$\alpha \in (\frac{\pi}{2},\pi), \hspace{0.5mm} V \geq 0$};
%\node at (-5.1,1.8) {\tiny{$-\frac{\pi}{2} < \alpha < 0 , \hspace{0.5mm} J = sign(V) = -$}};

\node at (0.7,0.2) {\tiny{$\times$}};
\node at (1.7,0.2) {\tiny{$\times$}};
\node at (1.4,0.5) {\tiny{$\times$}};
\node at (1.5,0.9) {\tiny{$\times$}};

\node at (-0.7,0.2) {\tiny{$\times$}};
\node at (-1.7,0.2) {\tiny{$\times$}};
\node at (-1.4,0.5) {\tiny{$\times$}};
\node at (-0.7,0.9) {\tiny{$\times$}};
\node at (-0.7,1.8) {\tiny{$\times$}};

\node at (0.1,0.9) {\tiny{$\times$}};
\node at (0.25,1.6) {\tiny{$\times$}};

\end{tikzpicture}
\caption{\textbf{Graphic illustration of the 
localization of the complex eigenvalues near a Landau level:} 
In a domain $\Omega_{q,\nu}^+(r^2,r_0^2) \cap \Omega_q^+(r^2,r_0^2) \cap 
E_q^+(\alpha,\theta)$, the number of complex eigenvalues of 
$\chh := \chho + e^{i\alpha}V$ is bounded by $\mathcal{O}( \vert \ln r \vert)$, 
see Theorem \ref{theo3}. For $\theta$ small enough and $0 < \varepsilon \leq 
\varepsilon_0$ small enough, $\chh_\varepsilon := \chho + \varepsilon e^{i\alpha}V$ 
has no complex eigenvalues in $\Omega_q^+(r^2,r_0^2) \cap E_q^+(\alpha,\theta)$. 
They are localized around the semi-axis $z = 2bq + e^{i(2\alpha - \pi)} ]0,+\infty)$, 
see \textbf{(i)} and \textbf{(iii)} of Theorem \ref{theo4}.}
\end{center}
\end{figure}
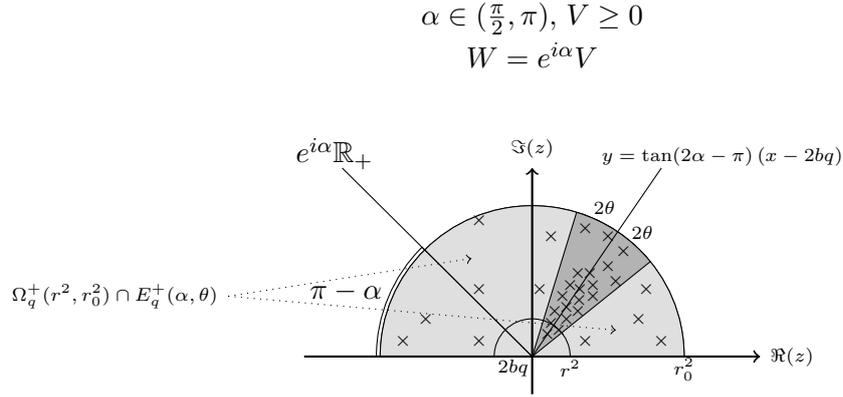

\noindent
Let us mention an important immediate consequence of Theorem 
\ref{theo4}-\textbf{(i)}.

\begin{cor}[Non-accumulation of complex eigenvalues]\label{c1}
Under the assumptions and the notations of Theorem \ref{theo4}, there is no 
accumulation of discrete eigenvalues of $\chh_\varepsilon$ near $2bq$, $q \in \bn$, for 
any $0 < \varepsilon \leq \varepsilon_0$, if
\begin{equation}
\alpha \in \pm \left( 0,\frac{\pi}{2} \right).
\end{equation}
\end{cor}

\noindent
Our results are summarized in Figure 2.2.
\\

\begin{figure}[h]\label{fig 1}
\begin{center}
\vspace*{-1cm}

\hspace*{-2.6cm} \includegraphics[scale=0.8]{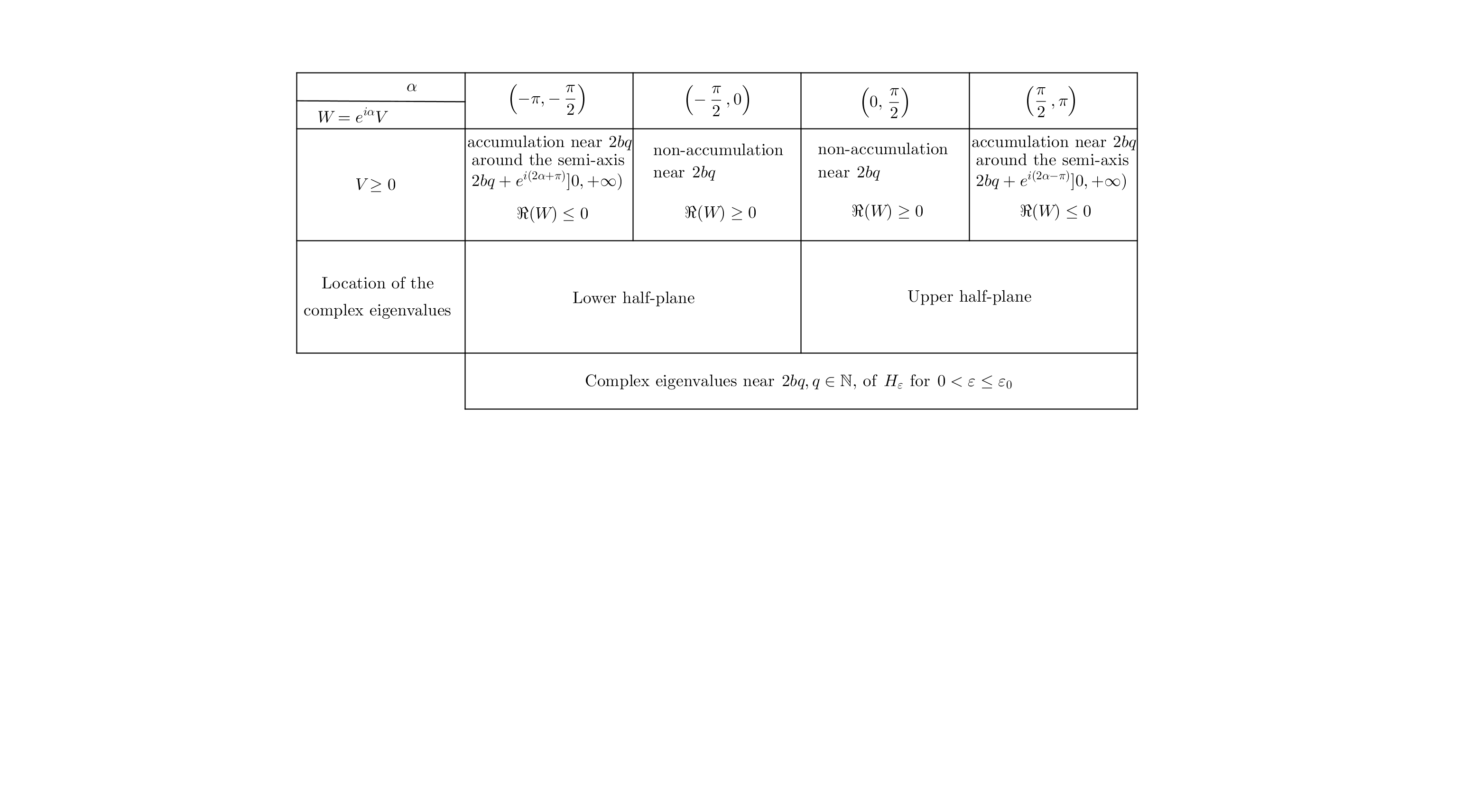}

\vspace*{-4.9cm}

\caption{Summary of results.}
\end{center}
\end{figure}

\noindent
About the accumulation of the complex eigenvalues of $\chh_\varepsilon$ near the landau 
levels, our results hold for each $0 < \varepsilon \leq \varepsilon_0$. Although this topic 
exceeds the scope of this paper, we expect this to be a general phenomenon in the sense of 
the following conjecture:

\begin{cnj}\label{cj2,1}
Let $W = \Phi V$ satisfy Assumption (A1) with $\Phi \in \bc \setminus \br e^{ik
\frac{\pi}{2}}$, $k \in \bz$, and $V : \br^3 \longrightarrow \br$ of definite 
sign. Then, there is no accumulation of complex eigenvalues of $\chh$ near $2bq$, 
$q \in\bn$, if and only if
$\Re(W) > 0$.
\end{cnj}

\noindent
The next result states that the Landau levels are the only possible accumulation points 
of the complex eigenvalues in some particular cases. Notations are those from above.

\begin{theo}[Dominated accumulation]\label{theo6}
Let the assumptions of Theorem \ref{theo4} hold with $\alpha \in \pm \left( 
\frac{\pi}{2},\pi \right)$. Then, for any $0 < \eta  < \sqrt{2b}$ and any $\theta > 0$ 
small enough, there exists $\Tilde{\varepsilon}_0 > 0$ such that for each $0 < \varepsilon 
\leq \Tilde{\varepsilon}_0$, $\chh_\varepsilon$ has no discrete eigenvalues in
\begin{equation}
\Omega_q^{\pm}(0,\eta^2) \setminus \left( 2bq + 
e^{i (2\alpha \mp \pi)} e^{i (-2\theta,2\theta)} (0,\eta^2) \right).
\end{equation}
If $\alpha \in \pm \left( 0,\frac{\pi}{2} \right)$, then $\chh_\varepsilon$ 
has no discrete eigenvalues in $\Omega_q^{\pm}(0,\eta^2)$. In particular, the Landau 
levels $2bq$, $q \in \bn$, are the only possible accumulation points of the 
discrete eigenvalues of $\chh_\varepsilon$.
\end{theo}

\begin{rem}
\textup{Since the Landau levels are the only possible accumulation points of the discrete
eigenvalues, then an immediate consequence of Theorem \ref{theo6} is that for $\alpha \in 
\pm \left( 0,\frac{\pi}{2} \right)$ there is no accumulation of complex eigenvalues of
$\chh_\varepsilon$, $0 < \varepsilon \leq \Tilde{\varepsilon}_0$, near the whole real 
axis.}
\end{rem}

\begin{rem}
\textup{In higher dimension $n \geq 3$, the magnetic self-adjoint Schrödinger operator 
$\chho$ in $L^{2}(\br^{n})$ has the form $(-i\nabla - \textbf{A})^{2}$, $\textbf{A} := 
(A_{1},\ldots,A_{n})$ being a magnetic potential generating the magnetic field. By 
introducing the $1$-form $\mathcal{A} := \sum_{j=1}^{n} A_{j}dx_{j}$, the magnetic field $\mathcal{\textbf{B}}$ can be defined as its exterior differential. Namely,
$\mathcal{\textbf{B}} := d \mathcal{A} = \sum_{j<k} B_{jk} dx_{j} 
\wedge dx_{k}$ with
\begin{equation}
\small{B_{jk} := 
\frac{\partial A_{k}}{\partial x_{j}} - 
\frac{\partial A_{j}}{\partial x_{k}}, \quad j,k = 1,\ldots,n.}
\end{equation}
For $n = 3$, the magnetic field is identified with $\textbf{B} = (B_{1},B_{2},B_{3}) := 
\text{curl} \hspace{0.1cm} \textbf{A}$, where $B_{1} = B_{23}$, $B_{2} = B_{31}$ and 
$B_{3} = B_{12}$. In the case where the $B_{jk}$ do not depend on $x \in \mathbb{R}^{n}$, 
the magnetic field can be viewed as a real antisymmetric matrix $\textbf{B} := \big\lbrace 
B_{jk} \big\rbrace_{j,k=1}^{n}$. Assume that $\textbf{B} \neq 0$, put $2d := \text{rank} \hspace{0.6mm} \textbf{B}$ and $k := n - 2d = \dim \text{Ker} \hspace{0.6mm} \textbf{B}$. 
Introduce $b_{1} \geq \ldots \geq b_{d} > 0$ the real numbers such that the non-vanishing eigenvalues of $\textbf{B}$ coincide with $\pm ib_{j}$, $j = 1,\ldots,d$. Consequently, 
in appropriate Cartesian coordinates $(x_{1},y_{1},\ldots, x_{d},y_{d}) \in \mathbb{R}^{2d} 
= \text{Ran} \hspace{0.6mm} \textbf{B}$ and $\lambda = (\lambda_{1},\ldots,\lambda_{k}) 
\in \mathbb{R}^{k} = \text{Ker} \hspace{0.6mm} \textbf{B}$, $k \geq 1$, the operator 
$\chho$ can be defined as
\begin{equation}\label{opg}
\small{\chho = \displaystyle \sum_{j=1}^{d} \left\lbrace \left( 
-i\frac{\partial}{\partial x_{j}} + \frac{b_{j}y_{j}}{2} \right)^{2} + 
\left( -i\frac{\partial}{\partial y_{j}} - \frac{b_{j}x_{j}}{2} \right)^{2} 
\right\rbrace + \sum_{\ell=1}^{k} \frac{\partial^{2}}{\partial 
\lambda_{\ell}^{2}}}.
\end{equation}
The operator $\chho$ given by \eqref{op0} considered in this paper, is just the magnetic 
Schrödinger operator defined by \eqref{opg} shifted by $-b$ in the particular case $n = 3$ 
(then $d = 1$, $k = 1$), $b_{1} = b_{2} = b$ and $b_3 = 0$. However, our results remain 
valid at least for the case $n = 2d + 1$ (then $k = 1$) with $d \geq 1$. The general case 
for the operator \eqref{opg} is an open problem.}
\end{rem}

\section{Preliminaries and first reductions}\label{s3}

\subsection{Schatten-von Neumann classes and determinants}\label{ss3.1}

\noindent
Recall that $\mathscr{H}$ denotes a separable Hilbert space. Let $\sinf(\mathscr{H})$ 
be the set of compact linear operators on $\mathscr{H}$. Denote $s_k(T)$ the $k$-th 
singular value of $T \in \sinf(\mathscr{H})$. The Schatten-von Neumann classes 
$\spp(\mathscr{H})$, $p \in [1,+\infty)$, are defined by
\begin{equation}\label{eq03,1}
\spp(\mathscr{H}) := \Big\lbrace T \in \sinf(\mathscr{H}) : 
\Vert T \Vert^p_\spp := \sum_k s_k(T)^p < +\infty \Big\rbrace.
\end{equation}
We will write simply $\spp$ when no confusion can arise. For $T \in \spp$, the 
$p$-regularized determinant is defined by
\begin{equation}\label{eq03,2}
\small{\textup{det}_{\lceil p \rceil} (I - T)
 := \prod_{\mu \hspace*{0.1cm} \in \hspace*{0.1cm} \sigma (T)} 
 \left[ (1 - \mu) \exp \left( \sum_{k=1}^{\lceil p \rceil-1} 
 \frac{\mu^{k}}{k} \right) \right]},
\end{equation}
where $\lceil p \rceil := \min \big\lbrace n \in \mathbb{N} : n \geq p \big\rbrace$. The 
following properties are well-known about this determinant \big(see for instance 
\cite{simo}\big):

\textbf{a)} $\textup{det}_{\lceil p \rceil} (I) = 1$.

\textbf{b)} For any bounded operators $A$, $B$ on $\mathscr{H}$ such
that $AB$ and $BA \in \spp$, $\textup{det}_{\lceil p \rceil} (I - AB)
= \textup{det}_{\lceil p \rceil} (I - BA)$.

\textbf{c)} The operator $I - T$ is invertible if and only if
$\textup{det}_{\lceil p \rceil} (I - T) \neq 0$.

\textbf{d)} If $T : \Omega \longrightarrow \spp$ is a holomorphic
operator-valued function on a domain $\Omega$, then so is the
function $\textup{det}_{\lceil p \rceil} \big( I - T(\cdot) \big)$ on $\Omega$.

\textbf{e)} If $T$ is a trace-class operator ($i.e.$ $T \in \mathcal{S}_1$), then 
\big(see for instance \cite[Theorem 6.2]{simo}\big)
\begin{equation}\label{eq03,3}
\small{\textup{det}_{\lceil p \rceil} (I - T)
= \textup{det} \hspace{0.05cm} (I - T) \exp \left( 
\sum_{k=1}^{\lceil p \rceil-1} \frac{\textup{Tr} \hspace{0.4mm} 
(T^{k})}{k} \right).}
\end{equation}

\textbf{f)} For $T \in \spp$, the inequality \big(see for instance 
\cite[Theorem 6.4]{simo}\big)
\begin{equation}\label{eq13,3}
\vert \textup{det}_{\lceil p \rceil} (I - T) \vert
\leq \exp \big( \Gamma_p \Vert T \Vert_\spp^p \big)
\end{equation}
holds, where $\Gamma_p$ is a positive constant depending only on $p$.

\textbf{g)} $\textup{det}_{\lceil p \rceil} (I - T)$ is Lipschitz as 
function on $\spp$ uniformly on balls:
\begin{equation}\label{eq23,3}
\big\vert \textup{det}_{\lceil p \rceil} (I - T_1) - 
\textup{det}_{\lceil p \rceil} (I - T_2) \big\vert
\leq \Vert T_1 - T_2 \Vert_\spp 
\exp \left( \Gamma_p \big( \Vert T_1 \Vert_\spp + \Vert T_2 
\Vert_\spp + 1 \big)^{\lceil p \rceil} \right),
\end{equation}
\big(see for instance \cite[Theorem 6.5]{simo}\big).

\subsection{On the relatively compactness of the potential $W$ with respect to $\chho$}

\begin{lem}\label{lem3,1}
Let $g \in L^p(\br^3)$, $p \geq 2$. Then, $g (\chho - z)^{-1} \in \spp$ for any 
$z \in \rho(\chho)$ (the resolvent set of $\chho$), with
\begin{equation}\label{eq3,0}
\small{\left\Vert g (\chho-z)^{-1} \right\Vert_{\spp}^{p} 
\leq  C \Vert g \Vert_{L^{p}}^{p} 
\textup{sup}_{s \in [0,+\infty)}^p \left\vert \frac{s + 1}{s - z} 
\right\vert},
\end{equation}
where $C = C(p)$ is constant depending on $p$.
\end{lem}

\noindent
\begin{proof}
Constants are generic, $i.e.$ changing from a relation to another.

First, let us show that \eqref{eq3,0} holds when $p$ is even. We have
\begin{equation}\label{eq3,1}
\left\Vert g (\chho-z)^{-1} \right\Vert_{\spp}^{p} 
\leq  \left\Vert g (\chho + 1)^{-1} \right\Vert_{\spp}^{p}
\left\Vert (\chho + 1)(\chho - z)^{-1} \right\Vert^p.
\end{equation}
By the Spectral mapping theorem,
\begin{equation}\label{eq3,2}
\left\Vert (\chho + 1)(\chho-z)^{-1} \right\Vert^p \leq 
\textup{sup}_{s \in [0,+\infty)}^p \left\vert \frac{s + 1}{s - z} 
\right\vert.
\end{equation}
Thanks to the resolvent identity, the diamagnetic inequality \big(see \cite[Theorem 2.3]{avr}-\cite[Theorem 2.13]{sim}\big) (only valid when $p$ is even), and the standard criterion \cite[Theorem 4.1]{sim}, 
we have
\begin{equation}\label{eq3,3}
\begin{split}
\left\Vert g \bigl( H_0 + 1 \bigr)^{-1} \right\Vert^p_\spp 
& \leq \left\Vert I + (H_0 + 1)^{-1}b \right\Vert^p
\left\Vert g \big( (-i\nabla - \textbf{A})^{2} + 1 \big)^{-1} 
\right\Vert^p_\spp \\
& \leq C \left\Vert g (-\Delta + 1)^{-1} \right\Vert^p_\spp
\leq C \Vert g \Vert_{L^{p}}^p 
\left\Vert \Bigl( \vert \cdot \vert^{2} + 1 \Bigr)^{-1} 
\right\Vert_{L^{p}}^{p}.
\end{split}
\end{equation}
So, for $p$ even, \eqref{eq3,0} follows by combining \eqref{eq3,1}, \eqref{eq3,2} and 
\eqref{eq3,3}.

Now, we show that \eqref{eq3,0} happens for any $p \geq 2$ by using interpolation method. 
If $p > 2$, there exists even integers $p_{0} < p_{1}$ such that $p \in (p_{0},p_{1})$ 
with $p_{0} \geq 2$. Let $s \in (0,1)$ satisfy $\frac{1}{p} = \frac{1-s}{p_{0}} + 
\frac{s}{p_{1}}$, and introduce the operator
$$
L^{p_i} \big(\mathbb{R}^{3}\big) \ni g \overset{T}{\longmapsto}  
g (H_{0} - z)^{-1} \in \mathcal{S}_{p_{i}}, \qquad i = 0, 1.
$$
Denote by $C_{i} = C(p_{i})$ the constant appearing in \eqref{eq3,0}, $i = 0$, $1$, and 
set 
$$
C(z,p_{i}) := C_i^{\frac{1}{p_i}} \textup{sup}_{s \in [0,+\infty)} \left\vert 
\frac{s + 1}{s - z} \right\vert.
$$
The inequality \eqref{eq3,0} implies that $\Vert T \Vert \leq C(z,p_{i})$ for $i = 0$, $1$. 
Now, with the help of the Riesz-Thorin Theorem \big(see for instance \cite[Sub. 5 of Chap. 
6]{fol}, \cite{rie,tho}, \cite[Chap. 2]{lun}\big), we can interpolate between $p_{0}$ and 
$p_{1}$ to get the extension 
$T : L^{p}\big( \br^{3}\big) \longrightarrow S_{p}$ with
$$
\Vert T \Vert \leq C(z,p_{0})^{1-s} C(z,p_{1})^{s} \leq 
C(p)^{\frac{1}{p}} \textup{sup}_{s \in [0,+\infty)} 
\left\vert \frac{s + 1}{s - z} \right\vert.
$$
In particular, for any $g \in L^p(\br^3)$, we have
$$
\Vert T(g) \Vert_{\spp} \leq 
C(p)^{\frac{1}{p}} \textup{sup}_{s \in [0,+\infty)} \left\vert 
\frac{s + 1}{s - z} \right\vert \Vert g \Vert_{L^p},
$$
which is equivalent to \eqref{eq3,0}. This concludes the proof.
\end{proof}

\noindent
Lemma \ref{lem3,1} above applied to the non-self-adjoint electric potential $W$ satisfying \textit{Assumption (A1)} gives
\begin{equation}\label{eq3,4}
\small{\left\Vert \vert W \vert^{\frac{1}{2}} (\chho-z)^{-1} \right\Vert_{\spp}^{p} 
\leq  C \Vert F \Vert_{L^{\frac{p}{2}}}^{\frac{p}{2}} \Vert G \Vert_{L^{\frac{p}{2}}}^{\frac{p}{2}}
\textup{sup}_{s \in [0,+\infty)}^p \left\vert \frac{s + 1}{s - z} \right\vert,
}
\end{equation}
for $p \geq 2$. In particular, $W$ is a relatively compact perturbation with respect $\chho$ 
since it is bounded.

\subsection{Reduction to zeros of a holomorphic function problem}

\noindent
Throughout this article, we deal with the following choice of the complex square root:
\begin{equation}\label{eq2,21}
\bc \setminus (-\infty,0] \overset{\sqrt{\cdot}}{\longrightarrow} 
\bc_+.
\end{equation}
For a fixed Landau level $\Lambda_{q} = 2bq$, $q \in \bn$, and $0 < \eta < \sqrt{2b}$, 
let $\Omega_q^\pm(0,\eta^2)$ be the half-rings defined by \eqref{eq2,3b}.
Make the change of variables $z - \Lambda_{q} = k^2$ and introduce
\begin{equation}\label{eq2,4}
\mathcal{D}_\pm^\ast (\eta) 
:= \big\lbrace k \in \bc_\pm : 0 < \vert k \vert < \eta : \Re(k) 
> 0 \big\rbrace.
\end{equation}

\begin{rem}\label{rc}
\textup{Notice that $\Omega_q^\pm(0,\eta^2)$ can be parametrized by 
$z = z_{q}(k) := \Lambda_{q} + k^2$ with $k \in \mathcal{D}_\pm^\ast (\eta)$ respectively 
(see Figure 3.1).}
\end{rem}

\begin{figure}[h]\label{fig 1}
\begin{center}

\vspace*{-1.5cm}

\hspace*{-1cm} \includegraphics[scale=0.6]{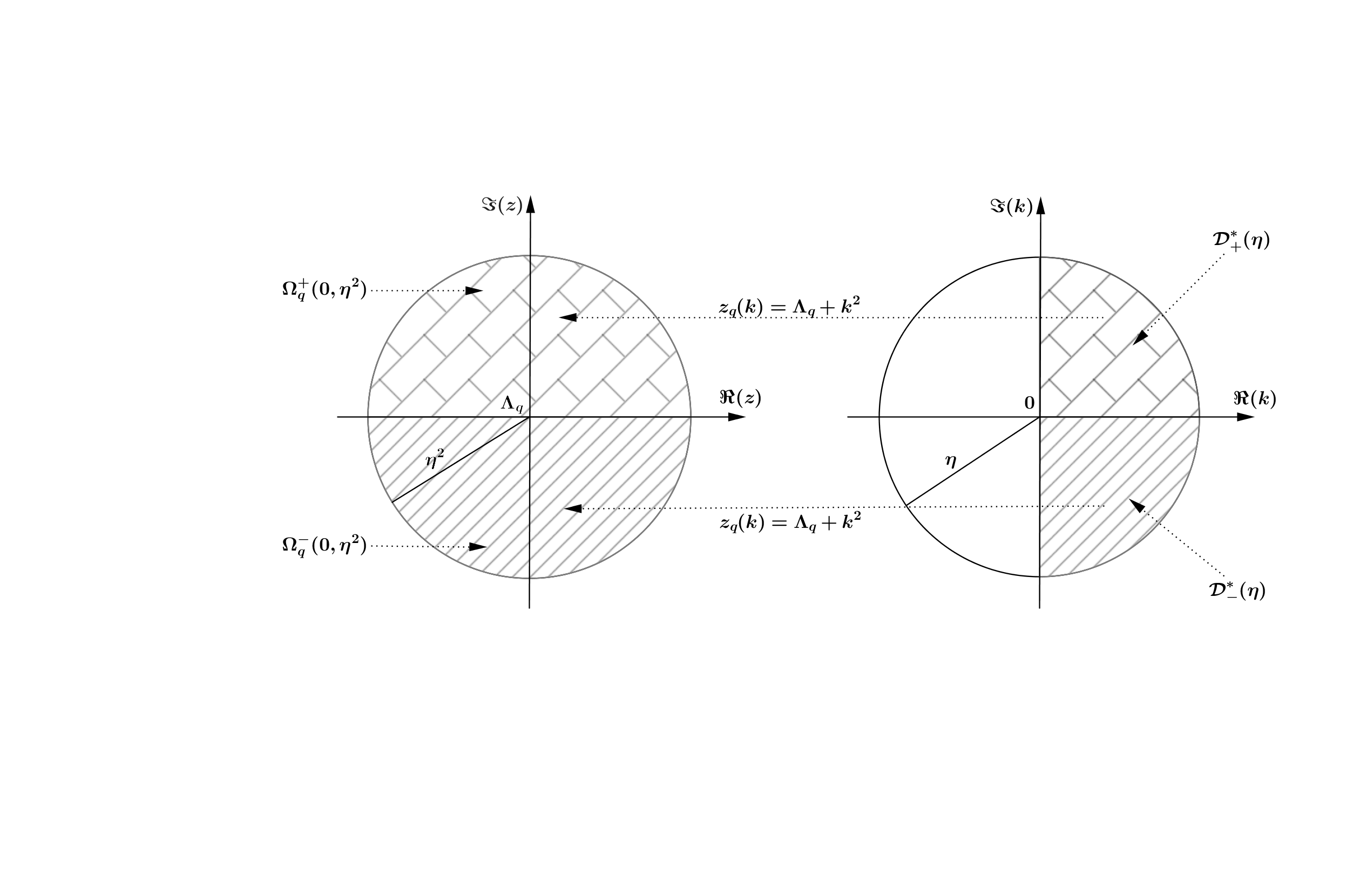}

\vspace*{-2cm}

\caption{Images $\Omega_q^\pm(0,\eta^2)$ of $\mathcal{D}_\pm^\ast(\eta)$ by the local parametrisation $z_q(k) = \Lambda_q + k^2$.}
\end{center}
\end{figure}

In this subsection, we show how we can reduce the investigation of the discrete 
eigenvalues $z_q(k) := \Lambda_q + k^2 \in \Omega_q^\pm(0,\eta^2)$, $k \in \mathcal{D}_\pm^\ast(\eta)$, to that of the zeros of a holomorphic function on 
$\Omega_q^\pm(0,\eta^2)$.

Let us recall that $P_q$, $q \in \bn$, is the projection onto 
${\rm Ker}(H_{\text{Landau}} - \Lambda_q)$. Hence, introduce in $L^2(\br^3)$ the 
projection $p_q := P_q \otimes I_3$, $q \in \bn$. With respect to the polar 
decomposition of $W$, write $W = \Tilde{J} \vert W \vert$. Then, for any 
$z \in \bc \setminus [0,+\infty)$, we have
\begin{equation}\label{eq3,5}
\begin{split}
& \Tilde{J} \vert W \vert^{\frac{1}{2}} (\chho - z)^{-1} \vert W 
\vert^{\frac{1}{2}} \\
& = \Tilde{J} \vert W \vert^{\frac{1}{2}} p_{q}(\chho - z)^{-1}\vert W 
\vert^{\frac{1}{2}}
+ \sum_{j \neq q} \Tilde{J} \vert W \vert^{\frac{1}{2}} 
p_{j}(\chho - z)^{-1} 
\vert W \vert^{\frac{1}{2}}.
\end{split}
\end{equation}
Since
$$
(H_0 - z)^{-1} = \sum_{q \in \mathbb{N}} P_q \otimes (D_{x_3}^2 + 
\Lambda_q - z)^{-1}, \qquad D_{x_3}^2 := -\frac{\partial^2}{\partial x_{3}^2},
$$
then for $z = z_q(k)$, $k \in \mathcal{D}_\pm^\ast(\eta)$, the identity 
\eqref{eq3,5} becomes
\begin{equation}\label{eq3,6}
\begin{split}
& \Tilde{J} \vert W \vert^{\frac{1}{2}} \big( \chho - z_q(k) 
\big)^{-1} \vert W \vert^{\frac{1}{2}} \\
& = \Tilde{J} \vert W \vert^{\frac{1}{2}} p_{q}(D_{x_3}^2 - 
k^2)^{-1}\vert W \vert^{\frac{1}{2}}
+ \sum_{j \neq q} \Tilde{J} \vert W \vert^{\frac{1}{2}} 
p_{j}(D_{x_3}^2 + \Lambda_j - \Lambda_q - k^2)^{-1} \vert 
W \vert^{\frac{1}{2}}.
\end{split}
\end{equation}
Hence, thanks to Lemma \ref{lem3,1}, we have the following 

\begin{prop}\label{prop3,1} 
Suppose that Assumption (A1) holds. Then, the operator-valued functions
$$
\mathcal{D}_\pm^\ast(\eta) \ni k \longmapsto \mathcal{T}_{W} \big( 
z_{q}(k) \big) := \Tilde{J} \vert W \vert^{\frac{1}{2}} \big( \chho - 
z_{q}(k) \big)^{-1} \vert W \vert^{\frac{1}{2}}
$$ 
are analytic with values in $\spp$.
\end{prop} 

\noindent
For $z \in \bc \setminus [0,+\infty)$, on account of Lemma \ref{lem3,1} and Subsection 
\ref{ss3.1}, we can introduce the $p$-regularized determinant 
$
\textup{det}_{\lceil p \rceil} \left( I + W(\chho - z)^{-1} 
\right).
$
The following characterization on the discrete eigenvalues is well known $\big($see for 
instance \cite[Chap. 9]{sim}$\big)$:
\begin{equation}\label{eq3,8}
z \in \sigma_{\textup{\textbf{disc}}}(\chh) \Leftrightarrow 
f(z) := \textup{det}_{\lceil p \rceil} \big( I + W(\chho - z)^{-1} \big) = 0,
\end{equation}
$\chh$ being the perturbed operator defined by \eqref{eq1,3}. According to Property 
\textbf{d)} of Subsection \ref{ss3.1}, if $W(\chho - \cdot)^{-1}$ is holomorphic on a 
domain $\Omega$, then so is the function $f$ on $\Omega$. Moreover, the algebraic  
multiplicity of $z \in \sigma_{\textup{\textbf{disc}}}(\chh)$ is equal to its order 
as zero of the function $f$.
\\

\noindent
In the next lemma, the notation $Ind_{\gamma}(\cdot)$ in the right hand-side of 
\eqref{eq3,9} is recalled in the Appendix.

\begin{prop}\label{prop3,2} 
Let $\mathcal{T}_{W} \big( z_{q}(k) \big)$ be defined by Proposition \ref{prop3,1}, 
$k \in \mathcal{D}_\pm^\ast(\eta)$. Then, the following assertions are equivalent:

$\textup{\textbf{(i)}}$ $z_{q}(k_{0}) := \Lambda_{q} + k_{0}^2 
\in \Omega_q^\pm(0,\eta^2)$ is a discrete eigenvalue of $\chh$,

$\textup{\textbf{(ii)}}$ 
$\textup{det}_{\lceil p \rceil} \big( I + \mathcal{T}_{W} \big( 
z_{q}(k_{0}) \big) \big) = 0$,

$\textup{\textbf{(iii)}}$ $-1$ is an eigenvalue of
$\mathcal{T}_{W} \big( z_{q}(k_{0}) \big)$.
\\
Moreover,
\begin{equation}\label{eq3,9}
\textup{mult} \big( z_{q}(k_{0}) \big) = 
Ind_{\gamma} \hspace{0.5mm} \Big( I + \mathcal{T}_{W}\big( 
z_{q}(\cdot) \big) \Big),
\end{equation}
$\gamma$ being a small contour positively oriented, containing $k_{0}$ as the unique 
point $k \in \mathcal{D}_\pm^\ast(\eta)$ verifying $z_{q}(k) \in \Omega_q^\pm(0,\eta^2)$ 
is a discrete eigenvalue of $\chh$.
\end{prop}

\noindent
\begin{proof}
The equivalence \textbf{(i)} $\Leftrightarrow$ \textbf{(ii)} is an immediate consequence 
of the characterization \eqref{eq3,8}, and the equality 
$$
\textup{det}_{\lceil p \rceil} \big( I + W(\chho - z)^{-1} \big)
= \textup{det}_{\lceil p \rceil} \big( I + \Tilde{J} \vert W 
\vert^{\frac{1}{2}} (\chho - z)^{-1} \vert W \vert^{\frac{1}{2}} 
\big).
$$

The equivalence \textbf{(ii)} $\Leftrightarrow$ \textbf{(iii)} is an obvious consequence of 
Property \textbf{c)} of Subsection \ref{ss3.1}.

Now, let us prove the equality \eqref{eq3,9}. Consider $f$ the function introduced in 
\eqref{eq3,8}. Thanks to the discussion just after \eqref{eq3,8}, if $\gamma'$ is a small 
contour positively oriented containing $z_{q}(k_{0})$ as the unique discrete eigenvalue 
of $\chh$, then
\begin{equation}\label{eq3,10}
\textup{mult} \big( z_{q}(k_{0}) \big) = ind_{\gamma'} f.
\end{equation}
The right hand-side of \eqref{eq3,10} being the index defined by \eqref{eqa,5} of the 
holomorphic function $f$ with respect to the contour $\gamma'$. Then, the equality 
\eqref{eq3,9} follows directly from the equality
$$
ind_{\gamma'} f = Ind_{\gamma} \hspace{0.5mm} 
\Big( I + \mathcal{T}_{W}\big( z_{q}(\cdot) \big) \Big),
$$
see for instance \cite[(2.6)]{bo} for more details.
This completes the proof.
\end{proof}

\section{Decomposition of the sandwiched resolvent}\label{s4}

We decompose the operator $\mathcal{T}_{W} \big( z_q(k) \big)$, $k \in 
\mathcal{D}_\pm^\ast(\eta)$, into a singular part at zero (corresponding to the singularity 
at the Landau level $\Lambda_q = 2bq$), and a holomorphic part in $\mathcal{D}_\pm^\ast(\eta)$, continuous on $\overline{\mathcal{D}_\pm^\ast(\eta)}$ with values in $\spp$.

First, note that due to our choice of the complex square root \eqref{eq2,21}, we 
respectively have $\sqrt{k^2} = \pm k$ for $k \in \mathcal{D}_\pm^\ast(\eta)$.

By \eqref{eq3,6}, we have
\begin{equation}\label{eq4,1}
\mathcal{T}_{W} \big( z_{q}(k) \big) = \Tilde{J} \vert W 
\vert^{\frac{1}{2}} 
p_{q}(D_{x_3}^2 - k^2)^{-1}\vert W \vert^{\frac{1}{2}}
+ \sum_{j \neq q} \Tilde{J} \vert W \vert^{\frac{1}{2}} 
p_{j}(D_{x_3}^2 + \Lambda_j - \Lambda_q - k^2)^{-1} 
\vert W \vert^{\frac{1}{2}}.
\end{equation}
Introduce $G_\pm$ the multiplication operators by the functions $G^{\pm \frac{1}{2}}(\cdot)$ respectively. Then, we have
\begin{equation}\label{eq4,2}
\Tilde{J} \vert W \vert^{\frac{1}{2}} 
p_{q}(D_{x_3}^2 - k^2)^{-1}\vert W \vert^{\frac{1}{2}}
= \Tilde{J} \vert W \vert^{\frac{1}{2}} G_- P_q \otimes G_+ 
(D_{x_3}^2 - k^2)^{-1}G_+ G_- \vert W \vert^{\frac{1}{2}}.
\end{equation}
It follows from the integral kernel 
\begin{equation}\label{ik} 
I_z(x_3,x_3') := -\frac{e^{i\sqrt{z} \vert x_3 - x_3' 
\vert}}{2i\sqrt{z}}
\end{equation}
of $(D_{x_3}^2 - z)^{-1}$, $\Im \big( \sqrt{z} \big) > 0$, that 
$G_+ (D_{x_3}^2 - k^2)^{-1}G_+$ admits the integral kernel 
\begin{equation}\label{eq4,3}
\pm G^{\frac{1}{2}}(x_3) \frac{i \textup{e}^{\pm i k \vert x_3 - 
x_3' \vert}}{2 k} G^{\frac{1}{2}}(x_3'), \quad k \in 
\mathcal{D}_\pm^\ast(\eta).
\end{equation}
Then, $G_+ (D_{x_3}^2 - k^2)^{-1}G_+$ can be decomposed as
\begin{equation}\label{eq4,4}
G_+ (D_{x_3}^2 - k^2)^{-1}G_+ = \pm \frac{1}{k}a + b(k), 
\quad k \in \mathcal{D}_\pm^\ast(\eta),
\end{equation}
where $a : L^{2}(\mathbb{R}) \longrightarrow L^{2}(\mathbb{R})$ is 
the rank-one operator defined by 
\begin{equation}\label{eq4,5}
a(u) := \frac{i}{2} \big\langle u,G_+ \big\rangle G_+,
\end{equation}
and $b(k)$ the operator with integral kernel
\begin{equation}\label{eq4,6}
\pm G^{\frac{1}{2}}(x_3) i \frac{ \textup{e}^{\pm i k \vert x_3 - x_3' 
\vert} - 1}{2 k} G^{\frac{1}{2}}(x_3').
\end{equation}
It can be easily remarked that $-2ia = c^\ast c$, where $c : L^{2}(\mathbb{R}) 
\longrightarrow \mathbb{C}$ is the operator defined by $c(u) := \langle u,G_+ \rangle$, 
so that $c^{\ast} : \mathbb{C} \longrightarrow L^{2}(\mathbb{R})$ is given 
by $c^{\ast}(\lambda) = \lambda G_+$. This together with \eqref{eq4,4}, \eqref{eq4,5}, 
and \eqref{eq4,6}, give for any $q \in \bn$,
\begin{equation}\label{eq4,7}
\begin{aligned}
P_q \otimes G_+ (D_{x_3}^2 - k^2)^{-1}G_+ = \pm \frac{i}{2k} P_q 
\otimes c^\ast c + P_q \otimes s(k), \quad k \in 
\mathcal{D}_\pm^\ast(\eta),
\end{aligned}
\end{equation}
where $s(k)$ is the operator acting from $G^{\frac{1}{2}}(x_3) 
L^{2}(\mathbb{R})$ to $G^{-\frac{1}{2}}(x_3) L^{2}(\mathbb{R})$ with integral kernel 
given by
\begin{equation}\label{eq4,8}
\pm \frac{ 1 - \textup{e}^{\pm i k \vert x_3 - x_3' \vert}}{2 i k}.
\end{equation}
By combining \eqref{eq4,2} and \eqref{eq4,7}, we get for any $k \in 
\mathcal{D}_\pm^\ast(\eta)$
\begin{equation}\label{eq4,9}
\begin{split}
& \Tilde{J} \vert W \vert^{\frac{1}{2}} 
p_{q}(D_{x_3}^2 - k^2)^{-1}\vert W \vert^{\frac{1}{2}} \\
& = \pm \frac{i\Tilde{J}}{2k} \vert W \vert^{\frac{1}{2}} G_- (P_q \otimes c^\ast c) 
G_- \vert W \vert^{\frac{1}{2}} 
+ \Tilde{J} \vert W \vert^{\frac{1}{2}} G_- P_q \otimes s(k) G_- \vert W \vert^{\frac{1}{2}}.
%& = \pm \frac{iJ}{2k} \left( (P_q \otimes c) G_{-} \vert W \vert^{\frac{1}{2}} \right)^{\ast} 
%\left( (P_q \otimes c) G_{-} \vert W \vert^{\frac{1}{2}} \right)
\end{split}
\end{equation}
The operator $P_q$ admits an explicit integral kernel
\begin{equation}\label{ker1}
\small{ {\mathcal P}_{q,b}(X_\perp,X_\perp^\prime)=\frac{b}{2\pi} 
L_q \left( \frac{b \vert X_\perp - X_\perp^\prime \vert^2}{2} \right)
\exp \Big( -\frac{b}{4} \big( \vert X_\perp - X_\perp^\prime \vert^2 + 
2i(x_1x_2^\prime - x_1^\prime x_2) \big) \Big) },
\end{equation}
where $X_\perp = (x_{1},x_{2}), X_\perp^\prime = (x_1^\prime,x_2^\prime) 
\in \br^2$, and $L_q (t) : = \frac{1}{q !} e^t \frac{d^q ( t^q e^{-t} )}{d t^q}$ are 
the Laguerre polynomials. Then, \eqref{eq4,9} becomes for any 
$k \in \mathcal{D}_\pm^\ast(\eta)$
\begin{equation}\label{eq4,10}
\Tilde{J} \vert W \vert^{\frac{1}{2}} 
p_{q}(D_{x_3}^2 - k^2)^{-1}\vert W \vert^{\frac{1}{2}}  = \pm 
\frac{i\Tilde{J}}{k} K^\ast K + \Tilde{J} \vert W \vert^{\frac{1}{2}} 
G_- P_q \otimes s(k) G_- \vert W \vert^{\frac{1}{2}},
\end{equation}
where the operator $K$ is given by
\begin{equation}\label{eq4,11}
K := \frac{1}{\sqrt{2}} (P_q \otimes c) G_{-} \vert W \vert^{\frac{1}{2}}.
\end{equation} 
To be more explicit, the operator $K : L^{2}(\br^{3}) \longrightarrow L^{2}(\br^{2})$ 
verifies
$$
(K \psi)(X_{\perp}) = \frac{1}{\sqrt{2}} 
\int_{\br^{3}} {\mathcal P}_{q,b}(X_{\perp},X_{\perp}^\prime) 
\vert W \vert^{\frac{1}{2}} (X_{\perp}^\prime,x_{3}^\prime) 
\psi (X_{\perp}^\prime,x_{3}^\prime) 
dX_{\perp}^\prime dx_{3}^\prime,
$$
${\mathcal P}_{q,b}(\cdot,\cdot)$ being the integral kernel given by \eqref{ker1}, 
while the adjoint operator $K^{\ast} : L^{2}(\br^{2}) \longrightarrow L^{2}(\br^{3})$ 
satisfies
$$
(K^{\ast}\varphi)(X_{\perp},x_{3}) = 
\frac{1}{\sqrt{2}} \vert W \vert^{\frac{1}{2}} (X_\perp,x_3) 
(P_q \varphi)(X_{\perp}).
$$
It is easy to check that 
$K K^{\ast} : L^{2}(\br^{2}) \longrightarrow L^{2}(\br^{2})$ verifies
\begin{equation}\label{eq4,12}
K K^{\ast} = P_{q} \textbf{\textup{W}} P_{q},
\end{equation}
where $\textbf{\text{W}}$ is the multiplication operator by the function 
$\textbf{\text{W}}$ defined by \eqref{eq2,1}. 

For $\lambda \in \br \setminus \lbrace 0 \rbrace$, define $(D_{x_3}^2 - \lambda)^{-1}$ 
as the operator with integral kernel 
\begin{equation}\label{eq4,121}
\displaystyle I_\lambda(x_3,x_3') := 
\lim_{\delta \downarrow 0} I_{\lambda + i\delta} (x_3,x_3') =
\left\{ \begin{array}{ccc} \frac{e^{-\sqrt{-\lambda} 
\vert x_3 - x_3' \vert}}{2\sqrt{-\lambda}} & \hbox{ if } & \lambda < 0, 
\\ -\frac{e^{i\sqrt{\lambda}\vert x_3 - x_3' \vert}}{2i\sqrt{\lambda}} & 
\hbox{ if } & \lambda > 0,  \end{array} \right.
\end{equation}
where $I_z(\cdot)$ is the integral kernel defined by \eqref{ik}. For $0 \leq \vert 
\lambda \vert < \sqrt{2b}$, we have
\begin{equation}\label{eq4,1211}
\begin{split}
\Big\Vert \sum_{j \neq q} \Tilde{J} \vert W 
\vert^{\frac{1}{2}} & p_{j}(D_{x_3}^2 + \Lambda_j - \Lambda_q - 
\lambda^2)^{-1} \vert W \vert^{\frac{1}{2}} \Big\Vert_\sd \\
\leq \sum_{j < q} & \Big\Vert \Tilde{J} \vert W 
\vert^{\frac{1}{2}} p_{j}(D_{x_3}^2 + \Lambda_j - \Lambda_q - 
\lambda^2)^{-1} \vert W \vert^{\frac{1}{2}} \Big\Vert_\sd \\
& + \Big\Vert \sum_{j > q} \Tilde{J} \vert W 
\vert^{\frac{1}{2}} p_{j}(D_{x_3}^2 + \Lambda_j - \Lambda_q - 
\lambda^2)^{-1} \vert W \vert^{\frac{1}{2}} \Big\Vert_\sd.
\end{split}
\end{equation}
Since $P_jP_\ell = \delta_{j,\ell}P_j$, then
\begin{equation}\label{eq4,122}
\begin{split}
\Big\Vert \sum_{j > q} \Tilde{J} & \vert W 
\vert^{\frac{1}{2}} p_{j}(D_{x_3}^2 + \Lambda_j - \Lambda_q - 
\lambda^2)^{-1} \vert W \vert^{\frac{1}{2}} \Big\Vert^2_\sd \\
& \leq \text{Const.} \sum_{j > q} \Big\Vert 
G_+ (D_{x_3}^2 + \Lambda_j - \Lambda_q - \lambda^2)^{-1} G_+
\Big\Vert^2_\sd.
\end{split}
\end{equation}
Using the integral kernel \eqref{eq4,121}, we obtain 
\begin{equation}\label{eq4,123}
\displaystyle
\left\{ \begin{array}{ccc} \left\Vert 
G_+ (D_{x_3}^2 + \Lambda_j - \Lambda_q - \lambda^2)^{-1}G_+ 
\right\Vert_\sd = \mathcal{O} \Big( \big\vert 2b(q - j) 
+ \lambda^2 \big\vert^{-\frac{1}{2}} \Big) & \hbox{ if } & j < q, \\  
\left\Vert G_+ (D_{x_3}^2 + \Lambda_j - \Lambda_q - \lambda^2)^{-1}
G_+ \right\Vert^2_\sd = \mathcal{O} \Big( \big\vert 2b(q - j) + \lambda^2 
\big\vert^{-\frac{3}{2}} \Big) & \hbox{ if } & j > q. 
\end{array} \right.
\end{equation}
This together with \eqref{eq4,122} imply that the left hand-side of \eqref{eq4,1211} 
is convergent in $\sd$. Moreover, arguing as in \cite[Proofs of Propositions 2.1-2.2]{fer},
it can be shown that 
$\overline{\mathcal{D}_\pm^\ast(\eta)} \ni k \mapsto 
\sum_{j \neq q} \Tilde{J} \vert W \vert^{\frac{1}{2}} p_{j}
(D_{x_3}^2 + \Lambda_j - \Lambda_q - k^2)^{-1} \vert W 
\vert^{\frac{1}{2}} \in \sd \big( L^{2}(\br) \big)$ 
is well defined and continuous. Similarly, as in \cite[Subsection 4.1]{bru}, it can be 
checked that $\overline{\mathcal{D}_\pm^\ast(\eta)} \ni k \mapsto 
G_+ s(k) G_+ \in \sd \big( L^{2}(\br) \big)$ is well defined and continuous. Therefore, 
the following proposition holds:

\begin{prop}\label{prop4,1} 
Assume that Assumption (A1) holds. Then, for $k \in \mathcal{D}_\pm^\ast(\eta)$,
\begin{equation}\label{eq4,13}
\mathcal{T}_{W} \big( z_{q}(k) \big) = \pm \frac{i\Tilde{J}}{k} 
\mathscr{B}_{q} + \mathscr{A}_{q}(k), \quad \mathscr{B}_{q} 
:= K^\ast K,
\end{equation}
where $\Tilde{J}$ is defined by the polar decomposition 
$W = \Tilde{J} \vert W \vert$. The operator 
$\mathscr{A}_{q}(k) \in \spp$ given by
\begin{align*}
\mathscr{A}_{q}(k) := 
\Tilde{J} \vert W \vert^{\frac{1}{2}} & G_- P_q \otimes s(k) G_- 
\vert W \vert^{\frac{1}{2}} \\ & + \sum_{j \neq q} \Tilde{J} 
\vert W \vert^{\frac{1}{2}} p_{j}(D_{x_3}^2 + \Lambda_j - 
\Lambda_q - k^2)^{-1} \vert W \vert^{\frac{1}{2}}
\end{align*}
is holomorphic on $\mathcal{D}_\pm^\ast(\eta)$ and continuous
on $\overline{\mathcal{D}_\pm^\ast(\eta)}$, $s(k)$ being 
defined by \eqref{eq4,7}.
\end{prop} 

\begin{rem}\label{r4.1}
\textup{\textbf{(i)} For any $r > 0$, we have}
\begin{equation}\label{eq4,14}
\textup{Tr} \hspace{0.4mm} \one_{(r,\infty)} 
\left( K^\ast K \right) 
= \textup{Tr} \hspace{0.4mm} \one_{(r,\infty)} 
\left( K K^\ast \right)
= \textup{Tr} \hspace{0.4mm} \one_{(r,\infty)} 
\big( P_{q} \textbf{\textup{W}} P_{q} \big).
\end{equation}

\textup{\textbf{(ii)} If $W$ satisfies \textit{Assumption (A2)} given by \eqref{eq2,5}, 
then Proposition \ref{prop4,1} holds with $\Tilde{J}$ replaced by $Je^{i\alpha}$, 
where $J = sign(V)$.}
\end{rem}

%\ref{theo2}
\section{Proof of Theorem \ref{theo2}: Upper bound, general case of electric potentials}\label{s5}

The proof falls into two parts.

\subsection{A preliminary Proposition}

\noindent
We begin by introducing the numerical range of $\chh$
$$
N(\chh) := \big\lbrace \langle \chh f,f \rangle : f \in Dom(\chh), 
\Vert f \Vert_{L^{2}} = 1 \big\rbrace.
$$ 
It is well known \big(see for instance \cite[Lemma 9.3.14]{dav}\big) that 
$\sigma(\chh) \subseteq \overline{N(\chh)}$.

\begin{prop}\label{prop4,2} 
Fix a Landau level $\Lambda_q := 2bq$, $q \in \bn$. Let $0 < s_{0} < \eta$ be 
sufficiently small. For any $k \in \lbrace 0 < s < \vert k \vert < s_{0} \rbrace \cap 
\mathcal{D}_\pm^\ast(\eta)$, 
 
$\textup{\textbf{(i)}}$ $z_{q}(k) := \Lambda_q + k^2$ is a discrete eigenvalue of $\chh$ 
near $\Lambda_{q}$ if and only if $k$ is a zero of
\begin{equation}\label{eq5,1}
\mathscr{D}(k,s) := \det \big( I + \mathscr{K}(k,s) \big).
\end{equation}
Here, $\mathscr{K}(k,s)$ is a finite-rank operator analytic with respect to $k$ verifying
$$
\small{\textup{rank} \hspace{0.6mm} \mathscr{K}(k,s) = \mathcal{O} 
\Big( \textup{Tr} \hspace{0.4mm} \one_{(s,\infty)} \big( P_q 
\textbf{\textup{W}} P_{q} \big) + 1 \Big),}
$$
and $\left\Vert \mathscr{K}(k,s) \right\Vert = \mathcal{O} \left( s^{-1} 
\right)$
uniformly with respect to $s < \vert k \vert < s_{0}$.

$\textup{\textbf{(ii)}}$ Further, if $z_{q}(k_{0}) := \Lambda_{q} + k_{0}^2$ is a discrete eigenvalue of $\chh$ near $\Lambda_q$, then
\begin{equation}\label{eq5,2}
\textup{mult} \big( z_{q}(k_{0}) \big) = Ind_{\gamma} \hspace{0.5mm} 
\left( I + \mathscr{K}(\cdot,s) \right) = \textup{m}(k_{0}),
\end{equation}
$\gamma$ being chosen as in \eqref{eq3,9} and $\textup{m}(k_{0})$ being the multiplicity 
of $k_{0}$ as zero of $\mathscr{D}(k,s)$.

$\textup{\textbf{(iii)}}$ If $z_{q}(k)$ verifies 
$\textup{dist} \big( z_q(k),\overline{N(\chh)} \big) > \varsigma > 0$, 
$\varsigma = \mathcal{O}(1)$, then $I + \mathscr{K}(k,s)$ is invertible and verifies
$
\left\Vert \big( I + \mathscr{K}(k,s) \big)^{-1} \right\Vert = 
\mathcal{O} \left( \varsigma^{-1} \right)
$
uniformly with respect to $s < \vert k \vert < s_{0}$.
\end{prop}

\noindent
\begin{proof}
\textbf{(i)-(ii)} Thanks to Proposition \ref{prop4,1}, $k \mapsto \mathscr{A}_{q}(k) 
\in \spp$ is continuous near zero. Thus, for $s_{0}$ sufficiently small, there exists a 
finite-rank operator $\mathscr{A}_{0}$ independent of $k$ and $\tilde{\mathscr{A}}(k) 
\in \spp$ continuous near zero, such that $\Vert \tilde{\mathscr{A}}(k) \Vert < \frac{1}{4}$, 
$\vert k \vert \leq s_{0}$, with
$$
\mathscr{A}_{q}(k) = \mathscr{A}_{0} + \tilde{\mathscr{A}}(k).
$$ 
Decompose $\mathscr{B}_{q}$ defined by \eqref{eq4,13} as
\begin{equation}\label{eq5,3}
\mathscr{B}_{q} = \mathscr{B}_{q} \one_{[0,\frac{1}{2}s]} (\mathscr{B}_{q}) 
+ \mathscr{B}_{q} \one_{]\frac{1}{2}s,\infty)} (\mathscr{B}_{q}).
\end{equation}
We have $\left\Vert \pm \frac{i\Tilde{J}}{k} \mathscr{B}_{q} 
\one_{[0,\frac{1}{2}s]} (\mathscr{B}_{q}) + \tilde{\mathscr{A}}(k) 
\right\Vert < \frac{3}{4}$ for $0 < s < \vert k \vert < s_{0}$ so that
\begin{equation}\label{eq5,4}
\small{\Big( I + \mathcal{T}_{W}\big( z_{q}(k) \big) \Big) = 
\big( I + \mathscr{K}(k,s) \big) \left( I \pm \frac{i\Tilde{J}}{k} 
\mathscr{B}_{q} \one_{[0,\frac{1}{2}s]} (\mathscr{B}_{q}) + 
\tilde{\mathscr{A}}(k) \right)},
\end{equation}
$\mathscr{K}(k,s)$ being given by
$$
\small{\mathscr{K}(k,s) := \left( \pm \frac{i\Tilde{J}}{k} 
\mathscr{B}_{q} \one_{]\frac{1}{2}s,\infty)} (\mathscr{B}_{q}) 
+ \mathscr{A}_{0} \right) \left( I \pm \frac{i\Tilde{J}}{k} 
\mathscr{B}_{q} \one_{[0,\frac{1}{2}s]} (\mathscr{B}_{q}) 
+ \tilde{\mathscr{A}}(k) \right)^{-1}}.
$$
Note that $\mathscr{K}(k,s)$ is a finite-rank operator. Moreover, 
thanks to \eqref{eq4,14}, its rank is of order
$$
\small{\mathcal{O} \Big( \textup{Tr} \hspace{0.4mm} 
\one_{(\frac{1}{2}s,\infty)} (\mathscr{B}_{q}) + 1 \Big) 
= \mathcal{O} \Big( \textup{Tr} \hspace{0.4mm} \one_{(s,\infty)} 
\big( P_q \textbf{\textup{W}} P_q \big) + 1 \Big)}.
$$
It is obvious that its norm is of order 
$\mathcal{O} \big( \vert k \vert^{-1} \big) = \mathcal{O} 
\big( s^{-1} \big)$. Since $\Vert \pm \frac{i\Tilde{J}}{k} 
\mathscr{B}_{q} \one_{[0,\frac{1}{2}s]} (\mathscr{B}_{q}) 
+ \tilde{\mathscr{A}}(k)\Vert < 1$ for $0 < s < \vert k 
\vert < s_{0}$, then
$$
\small{Ind_{\gamma} \hspace{0.5mm} \left(I \pm \frac{i\Tilde{J}}{k} 
\mathscr{B}_{q} \one_{[0,\frac{1}{2}s]} (\mathscr{B}_{q}) 
+ \tilde{\mathscr{A}}(k) \right) = 0}
$$ 
by \cite[Theorem 4.4.3]{goh}. Hence, \eqref{eq5,2} follows by applying to \eqref{eq5,4} 
the properties of the index of a finite meromorphic operator-valued function given in the 
Appendix. Thus, Proposition \ref{prop3,2} together with \eqref{eq5,2} show that $z_{q}(k)$ 
is a discrete eigenvalue of $\chh$ if and only if $k$ is a zero of the determinant 
$\mathscr{D}(k,s)$ defined by \eqref{eq5,1}.

\textbf{(iii)} Thanks to \eqref{eq5,4}, for $0 < s < \vert k \vert < s_{0}$, we have
\begin{equation}\label{eq5,5}
\small{I + \mathscr{K}(k,s) = \Big( I + \mathcal{T}_{W}\big( 
z_{q}(k) \big) \Big) \left( I + \frac{\Tilde{J}}{k} 
\mathscr{B}_{q} \one_{[0,\frac{1}{2}s]} (\mathscr{B}_{q}) 
+ \tilde{\mathscr{A}}(k) \right)^{-1}}.
\end{equation}
It is easy to check from the resolvent equation that
$$
\small{\left( I + \Tilde{J} \vert W \vert^{1/2} (\chho - z)^{-1} 
\vert W \vert^{1/2} \right) \left( I - \Tilde{J} \vert W 
\vert^{1/2} (\chh - z)^{-1} \vert W \vert^{1/2} \right) = I}.
$$
Thus, if $z_{q}(k)$ belongs to the resolvent set of $\chh$,
then
$$
\Big( I + \mathcal{T}_{W}\big( z_{q}(k) \big) \Big)^{-1} 
= I - \Tilde{J} \vert W \vert^{1/2} \big( \chh - z_{q}(k) 
\big)^{-1} \vert W \vert^{1/2}.
$$
Consequently, according to \eqref{eq5,5}, the operator $I + \mathscr{K}(k,s)$ is 
invertible for $0 < s < \vert k \vert < s_{0}$, and thanks to \cite[Lemma 9.3.14]{dav},
its satisfies for $\textup{dist} \big( z_q(k),\overline{N(\chh)} \big) > \varsigma > 0$,
$\varsigma = \mathcal{O}(1)$,
\begin{align*}
\left\Vert \big( I + \mathscr{K}(k,s) \big)^{-1} \right\Vert 
& = \mathcal{O} \Big( 1 + \left\Vert \vert W \vert^{1/2} 
\big( \chh - z_{q}(k) \big)^{-1} \vert W \vert^{1/2} \right\Vert \Big) \\
&= \mathcal{O} \big( 1 + \textup{dist} \big( z_q(k),\overline{N(\chh)} \big)^{-1} \big) \\
& = \mathcal{O} \left( \varsigma^{-1} \right).
\end{align*}
This concludes the proof.
\end{proof}

%\ref{theo2}
\subsection{Back to the proof of Theorem \ref{theo2}}

\noindent
Thanks to Proposition \ref{prop4,2}, for any $0 < s < \vert k \vert < s_{0}$,
\begin{equation}\label{eq5,6}
\begin{aligned}
\mathscr{D}(k,s) & = \prod_{j=1}^{\mathcal{O} \big( \textup{Tr} 
\hspace{0.4mm} \one_{(s,\infty)} (P_{q} \textbf{\textup{W}} 
P_{q}) + 1 \big)} \big{(} 1 + \lambda_{j}(k,s) \big{)} \\
& = \mathcal{O}(1) \hspace{0.5mm} \textup{exp} \hspace{0.5mm} 
\Big( \mathcal{O} \big( \textup{Tr} \hspace{0.4mm} 
\one_{(s,\infty)} \big( P_{q} \textbf{\textup{W}} P_{q} \big) 
+ 1 \big) \vert \ln s \vert \Big),
\end{aligned}
\end{equation}
where the $\lambda_{j}(k,s)$ are the eigenvalues of the operator $\mathscr{K} := 
\mathscr{K}(k,s)$ verifying $\vert \lambda_{j}(k,s) \vert = \mathcal{O} \left( s^{-1} 
\right)$. Consider $z_q(k)$ satisfying $\textup{dist} \big( z_q(k),\overline{N(\chh)} 
\big) > \varsigma > 0$ and $0 < s < \vert k \vert < s_{0}$. We have
$$
\mathscr{D}(k,s)^{-1} = \det \big( I + \mathscr{K} \big)^{-1} 
= \det \big( I - \mathscr{K} ( I + \mathscr{K})^{-1} \big),
$$
and as in \eqref{eq5,6}, we can show that
\begin{equation}\label{eq5,7}
\small{\vert \mathscr{D}(k,s) \vert \geq C \hspace{0.5mm} 
\textup{exp} \hspace{0.5mm} \Big( - C \big( \textup{Tr} 
\hspace{0.4mm} \one_{(s,\infty)} \big( P_{q} \textbf{\textup{W}} 
P_{q} \big) + 1 \big) \big( \vert \textup{ln} \hspace{0.5mm} 
\varsigma \vert + \vert \textup{ln} \hspace{0.5mm} s \vert \big) 
\Big)}.
\end{equation}
In particular, for $s^2 < \varsigma < 4s^2$, $0 < r \ll 1$, we deduce from 
\eqref{eq5,7} that
\begin{equation}\label{eq5,70}
- \ln \, \vert \mathscr{D}(k,s) \vert \leq
C \, \textup{Tr} \hspace{0.4mm} \one_{(s,\infty)} 
\big( P_{q} \textbf{\textup{W}} P_{q} \big) \vert \textup{ln} \hspace{0.5mm} s \vert
+ \mathcal{O}(1).
\end{equation}
Now, consider the domains $\Delta _{\pm} := k \in 
\big\lbrace r < \vert k \vert < 2r : \vert \Re(k) \vert > \sqrt{\frac{\nu}{2}} 
: \vert \Im(k) \vert > \sqrt{\frac{\nu}{2}} \big\rbrace \cap 
\mathcal{D}_{\pm }^{\ast }(\eta )$ with 
$0 < r < \sqrt{\Vert W \Vert_\infty} < \sqrt{\frac{5}{2}}r$ and 
$0 < \nu < 2r^2$. Since the numerical range $N(\chh)$ the operator 
$\chh$ is such that
\begin{equation}
N(\chh) \subseteq \big\lbrace z \in \bc : \vert \Im(z) \vert \leq 
\Vert W \Vert_\infty \big\rbrace,
\end{equation}
then we can find some $k_{0} \in\Delta _{\pm }/r$ such that $\textup{dist} 
\big( z_q(rk_0), \overline{N(\chh)} \big) \geq \varsigma > r^2$, $\varsigma < 4r^2$. 
Applying the Jensen Lemma \ref{la,1} with the function $g(k) := \mathscr{D}(rk,r)$, 
$k \in \Delta _{\pm}/r$, together with \eqref{eq5,6} and \eqref{eq5,70}, we get 
immediately Theorem \ref{theo2}.

\section{Proof of Theorem \ref{theo3}: Upper bound, special case of electric  potentials}\label{s6}

We prove only the case $\alpha = \frac{3\pi}{4}$; the case $\alpha = -\frac{3\pi}{4}$ 
follows similarly by replacing $k$ by $-k$, according to \textbf{(ii)} of Remark \ref{r2}
together with Remark \ref{rc} and \textbf{(ii)} of Remark \ref{r4.1}.
\\

\noindent
For any $\theta > 0$ small enough, set $\delta = \tan (\theta)$. Introduce the sector
\begin{equation}\label{eq2,6}
\mathcal{C}_\delta := \big\lbrace k \in \bc : - \delta \Im(k) \leq
\vert \Re(k) \vert \big\rbrace.
\end{equation}

\noindent
Let the assumptions of Theorem \ref{theo3} hold. Then, by Remark \ref{r4.1}, for any 
$q \in \bn$,
\begin{equation}\label{eq6,1}
\mathcal{T}_{W} \big( z_{q}(k) \big) = \frac{ie^{i\alpha}}{k} 
\mathscr{B}_{q}+ \mathscr{A}_{q}(k), \qquad
k \in \mathcal{D}_+^\ast(\eta),
\end{equation}
where $\mathscr{B}_{q}$ is a positive self-adjoint operator which does not depend 
on $k$, and $\mathscr{A}_{q}(k) \in \spp$ is continuous near $k = 0$. Denote 
$r_+ := \max (r,0)$. Since $I + \frac{ie^{i\alpha}}{k} \mathscr{B}_{q} = 
\frac{ie^{i\alpha}}{k} (\mathscr{B}_{q} - ike^{-i\alpha})$, then for $ike^{-i\alpha} 
\notin \sigma (\mathscr{B}_{q})$, the operator $I + \frac{ie^{i\alpha}}{k}\mathscr{B}_{q}$ 
is invertible with
\begin{equation}\label{eq6,2}
\small{\left\Vert \left( I + \frac{ie^{i\alpha}}{k} \mathscr{B}_{q} 
\right)^{-1} \right\Vert \leq \frac{\vert k \vert}{\sqrt{\big( 
\Im(ke^{-i\alpha}) \big)_+^2 + \vert \Re(ke^{-i\alpha}) \vert^2}}}.
\end{equation}
Moreover, for $k \in e^{i\alpha} \mathcal{C}_\delta$, it can be checked that,
uniformly with respect to $k$, $0 < \vert k \vert < r_0$,
\begin{equation}\label{eq6,3}
\small{\left\Vert \left( I + \frac{ie^{i\alpha}}{k} \mathscr{B}_{q} 
\right)^{-1} \right\Vert \leq \sqrt{1 + \delta^{-2}}}.
\end{equation}
Then, according to \eqref{eq6,1}, we can write
\begin{equation}\label{eq6,4}
\small{I + \mathcal{T}_{W} \big( z_{q}(k) \big) = 
\big( I + A(k) \big) \left( I + \frac{ie^{i\alpha}}{k} 
\mathscr{B}_{q} \right)},
\end{equation}
where 
\begin{equation}\label{eq6,5}
\small{A(k) := \mathscr{A}_{q} (k)
\left( I + \frac{ie^{i\alpha}}{k} \mathscr{B}_{q} \right)^{-1}} 
\in \spp.
\end{equation}
An easy computation shows that
$$
\mathcal{T}_{W} \big( z_{q}(k) \big) - A(k) = \big( I + A(k) 
\big) \frac{ie^{i\alpha}}{k} \mathscr{B}_{q} \in \mathcal{S}_1,
$$
since $\mathscr{B}_{q}$ is a trace-class operator if the function $F$ in 
\textit{Assumption (A1)} satisfies $F \in L^1(\br^2)$. Then, we get for
any $n \in \bn^\ast$
\begin{equation}\label{eq6,51}
\small{\mathcal{T}_{W}^n - A^n = \mathcal{T}_{W}^{n-1} \left( 
\mathcal{T}_{W} - A \right) + \left( \mathcal{T}_{W}^{n-1} - 
A^{n-1} \right) A} \in \mathcal{S}_1,
\end{equation}
So, by approximating $A(k)$ by a finite rank-operator and using the fact that
$$
\small{\textup{det}_{\lceil p \rceil} (I + T) = \textup{det} 
(I + T) \exp \left( \sum_{n=1}^{\lceil p \rceil-1} 
\frac{ (-1)^n \textup{Tr} \hspace{0.4mm} (T^n) }{n} \right)}
$$
for a trace-class operator $T$ \big(see Property \textbf{e)} of 
Subsection \ref{ss3.1} given by \eqref{eq03,3}\big), it can be shown with the help 
of \eqref{eq6,4} that
\begin{equation}\label{eq6,6}
\begin{split}
\small{\textup{det}_{\lceil p \rceil} \big( I +} & 
\small{\mathcal{T}_{W} \big( z_{q}(k) \big) \big)
= \textup{det} \left( I + \frac{ie^{i\alpha}}{k} \mathscr{B}_{q} 
\right)} \\
& \small{\times \textup{det}_{\lceil p \rceil} \big( I + A(k) 
\big) \exp \left( \sum_{n=1}^{\lceil p \rceil-1} \frac{ (-1)^n 
\textup{Tr} \hspace{0.4mm} \bigl( \mathcal{T}_{W}^n - A^n 
\bigr)}{n} \right)}.
\end{split}
\end{equation}
Thus, for $0 < \vert k \vert < r_0$ small enough, $k \in e^{i\alpha} 
\mathcal{C}_\delta$, the zeros of $\textup{det}_{\lceil p \rceil} \big( I + 
\mathcal{T}_{W} \big( z_{q}(k) \big) \big)$ are those of $\textup{det}_{\lceil p 
\rceil} \big( I + A(k) \big)$ with the same multiplicities thanks to Proposition 
\ref{prop3,2} and Property \eqref{eqa,52} applied to \eqref{eq6,4}.
\\

\noindent
Since $\mathscr{A}_{q}(\cdot) \in \spp$ is continuous near $k = 0$, this together with 
\eqref{eq6,3} implies that the $\spp$-norm of $A(k)$ is uniformly bounded with respect 
to $0 < \vert k \vert < r_0$ small enough, $k \in e^{i\alpha} \mathcal{C}_\delta$. Then, 
thanks to Property \textbf{f)} of Subsection \ref{ss3.1} given by \eqref{eq13,3}, we
have
\begin{equation}\label{eq6,7}
\small{\textup{det}_{\lceil p \rceil} \big( I + A(k) \big)
= \mathcal{O} \left( e^{\mathcal{O} 
\big( \Vert A(k) \Vert_{\spp}^p \big)} \right) = \mathcal{O}(1).}
\end{equation}

Now, let us establish a lower bound of $\textup{det}_{\lceil p \rceil} \big( I + A(k) 
\big)$. Thanks to \eqref{eq6,4}, we have
\begin{equation}\label{eq6,8}
\small{ \big( I + A(k) \big)^{-1} = \left( I + \frac{ie^{i\alpha}}{k} \mathscr{B}_{q} 
\right) \big( I + \mathcal{T}_{W} \big( z_{q}(k) \big)^{-1}}.
\end{equation}
Hence, by reasoning as in the proof of \textbf{(iii)}-Proposition \ref{prop4,2}, we 
obtain for $0 < s < \vert k \vert < r_0$ and $\textup{dist} \big( z_q(k),\overline{N(\chh)} 
\big) > \varsigma > 0$, $\varsigma = \mathcal{O}(1)$, uniformly with respect to $(k,s)$,
\begin{equation}\label{eq6,9}
\small{\left\Vert \big( I + A(k) \big)^{-1} \right\Vert
= \mathcal{O} \big( s^{-1} \big)
\mathcal{O} \big( \varsigma^{-1} \big)}.
\end{equation}
Let $(\mu_j)_j$ be the sequence of eigenvalues of $A(k)$. We have
\begin{equation}\label{eq6,10}
\begin{aligned}
\small{\left\vert \big( \textup{det}_{\lceil p \rceil} 
(I + A(k)) \big)^{-1} \right\vert}
& \small{= \left\vert \textup{det} \left( (I + A(k))^{-1} 
e^{\sum_{n=1}^{\lceil p \rceil-1} 
\frac{ (-1)^{n+1} A(k)^n}{n}} \right) \right\vert} \\
& \small{\leq \prod_{\vert \mu_j \vert \leq \frac{1}{2}}  
\left\vert \frac{e^{\sum_{n=1}^{\lceil p \rceil-1} 
\frac{ (-1)^{n+1} \mu_j^n}{n}}}{1 + \mu_j} \right\vert 
\times 
\prod_{\vert \mu_j \vert > \frac{1}{2}}  
\frac{e^{\left\vert \sum_{n=1}^{\lceil p \rceil-1} 
\frac{ (-1)^{n+1} \mu_j^n}{n} \right\vert}}{\vert 1 + \mu_j 
\vert}.
}
\end{aligned}
\end{equation}
Using the fact that $A(k)$ is uniformly bounded in $\spp$ with respect to 
$0 < \vert k \vert < r_0$ small enough, $k \in e^{i\alpha} \mathcal{C}_\delta$, 
it is easy to check that the first product is uniformly bounded. On the 
other hand, thanks to \eqref{eq6,9}, we have for $0 < s < \vert k \vert < r_0$ 
and $\textup{dist} \big( z_q(k),\overline{N(\chh)} \big) > \varsigma > 0$, 
$\varsigma = \mathcal{O}(1)$,
\begin{equation}\label{eq6,11}
\small{\vert 1 + \mu_j \vert^{-1}
= \mathcal{O} \big( s^{-1} \big)
\mathcal{O} \big( \varsigma^{-1} \big)},
\end{equation}
uniformly with respect to $k$, $s$. Consequently, since there exists a finite number 
of terms $\mu_j$ lying in the second product, we deduce from \eqref{eq6,10} that
\begin{equation}\label{eq6,12}
\small{\left\vert \textup{det}_{\lceil p \rceil} \big( I + A(k) \big) \right\vert 
\geq C e^{-C \big( \vert \ln \varsigma \vert + \vert \ln s \vert \big)}},
\end{equation}
for some positive constant $C > 0$. Now, one concludes as in the proof of Theorem 
\ref{theo2} by using the Jensen Lemma \ref{la,1}.

%\ref{theo4}, \ref{theo5}
%\section{Proof of Theorems $\ref{theo4}$ and $\ref{theo6}$}\label{s7}

\section{Theorem \ref{theo4}: Lower bound, upper bound and sectors free of complex eigenvalues}\label{s7}

As in the previous section, we only prove the case $\alpha \in (0,\pi)$. For 
$\alpha \in -(0,\pi)$, it suffices to replace $k$ by $-k$.

\textbf{(i)} Under the assumptions of Theorem \ref{theo4}, for any $q \in \bn$, 
we have
\begin{equation}\label{eq7,1}
I + \mathcal{T}_{\varepsilon W} \big( z_{q}(k) \big) = 
I + \frac{i\varepsilon e^{i\alpha}}{k} \mathscr{B}_{q} + 
\varepsilon \mathscr{A}_{q}(k), \qquad k \in \mathcal{D}_+^\ast(\eta).
\end{equation}
Similarly to the proof of Theorem \ref{theo3}, for $ike^{-i\alpha} \notin \sigma 
(\varepsilon \mathscr{B}_{q})$, the operator $I + \frac{i\varepsilon e^{i\alpha}}{k} \mathscr{B}_{q}$ is invertible. Further, for $k \in e^{i\alpha} \mathcal{C}_\delta$, $\delta = \tan (\theta)$,
we have
\begin{equation}\label{eq7,2}
\small{\left\Vert \left( I + \frac{i\varepsilon e^{i\alpha}}{k} \mathscr{B}_{q} 
\right)^{-1} \right\Vert \leq \sqrt{1 + \delta^{-2}}},
\end{equation}
uniformly with respect to $k$, $0 < \vert k \vert < r_0$. Then, as in \eqref{eq6,4} 
and \eqref{eq6,5}, we have
\begin{equation}\label{eq7,3}
\small{I + \mathcal{T}_{\varepsilon W} \big( z_{q}(k) \big) = 
\big( I + A(k) \big) \left( I + \frac{i\varepsilon e^{i\alpha}}{k} 
\mathscr{B}_{q} \right)},
\end{equation}
with
\begin{equation}\label{eq7,4}
\small{A(k) := \varepsilon \mathscr{A}_{q} (k) \left( I + \frac{i\varepsilon 
e^{i\alpha}}{k} \mathscr{B}_{q} \right)^{-1}} \in \spp.
\end{equation}
Since $\mathscr{A}_{q}(\cdot) \in \spp$ is continuous near $k = 0$, then there exists 
a constant $C > 0$ such that $\Vert \mathscr{A}_{q}(k) \Vert \leq C$. This together 
with \eqref{eq7,2} and \eqref{eq7,4} imply that for 
$0 < \varepsilon < \big( C \sqrt{1 + \delta^{-2}} \big)^{-1}$, the operator $I + \mathcal{T}_{\varepsilon W} \big( z_{q}(k) \big)$ is invertible for 
$k \in e^{i\alpha} \mathcal{C}_\delta$. Consequently, $z_q(k)$ is not a discrete 
eigenvalue.

\medskip

\textbf{(ii)} Decompose $\varepsilon\mathscr{B}_{q}$ as 
$\varepsilon\mathscr{B}_{q} = \mathscr{B}_+ + \mathscr{B}_-$, where 
$\mathscr{B}_+$ and $\mathscr{B}_-$ are defined by 
\begin{equation}\label{eq7,5}
\mathscr{B}_+ := \varepsilon\mathscr{B}_{q} \one_{[\frac{r}{2},4r]} 
(\varepsilon\mathscr{B}_{q}), \qquad 
\mathscr{B}_- := \varepsilon\mathscr{B}_{q} \one_{]0,\frac{r}{2}[ 
\cup ]4r,\infty[} (\varepsilon\mathscr{B}_{q}).
\end{equation}
It is easy to verify that for $\frac{2r}{3} < \vert k \vert < 
\frac{3r}{2}$, we have $\sigma \big( \frac{1}{\vert k \vert} 
\mathscr{B}_- \big) \subset \left[ 0,\frac{3}{4} \right] \cup 
\left[ \frac{8}{3},\infty \right[$. Therefore, $I + \frac{ie^{i\alpha}}{k}
\mathscr{B}_-$ is invertible with
\begin{equation}\label{eq7,6}
\small{\left\Vert \left( I + \frac{i e^{i\alpha}}{k} 
\mathscr{B}_- \right)^{-1} \right\Vert \leq 4},
\end{equation}
uniformly with respect to $0 < \vert k \vert < r_0$. Thus, for 
$0 < \varepsilon \leq \varepsilon_0$ small enough, $I + \frac{ie^{i\alpha}}{k} 
\mathscr{B}_- + \varepsilon \mathscr{A}_{q}(k)$ is invertible with a 
uniformly bounded inverse given by
\begin{equation}\label{eq7,7}
\small{\left( I + \frac{ie^{i\alpha}}{k} \mathscr{B}_- + 
\varepsilon \mathscr{A}_{q}(k) \right)^{-1} = \left( I + \frac{ie^{i\alpha}}{k} 
\mathscr{B}_- \right)^{-1} \left( I + \varepsilon \mathscr{A}_{q}(k) 
\left( I + \frac{ie^{i\alpha}}{k} \mathscr{B}_- \right)^{-1} \right)^{-1}.}
\end{equation}
This together with \eqref{eq7,1} and \eqref{eq7,5} allow to write
\begin{equation}\label{eq7,8}
\small{I + \mathcal{T}_{\varepsilon W} \big( z_{q}(k) \big) = 
\left( I + \frac{i e^{i\alpha}}{k} \mathscr{B}_- + \varepsilon 
\mathscr{A}_{q}(k) \right) \left( I + \left( I + \frac{i 
e^{i\alpha}}{k} \mathscr{B}_- + \varepsilon \mathscr{A}_{q}(k) 
\right)^{-1} \frac{i e^{i\alpha}}{k} \mathscr{B}_+ \right).
}
\end{equation}
Since $I + \frac{i e^{i\alpha}}{k} \mathscr{B}_- + \varepsilon \mathscr{A}_{q}(k)$ 
is invertible and $\mathscr{B}_+$ is a trace-class operator, then by exploiting Proposition 
\ref{prop3,2} and Property \eqref{eqa,52} applied to \eqref{eq7,8}, we see 
that the discrete eigenvalues of $\chh_\varepsilon$ are the zeros of 
\begin{equation}\label{eq7,9}
\small{\Tilde{D}(k,r) := 
\textup{det} \left( I + \left( I + \frac{i e^{i\alpha}}{k} \mathscr{B}_- + 
\varepsilon \mathscr{A}_{q}(k) \right)^{-1} \frac{i e^{i\alpha}}{k} \mathscr{B}_+ 
\right)}
\end{equation}
with the same multiplicities. Moreover, since $\frac{i e^{i\alpha}}{k} \mathscr{B}_+$ 
is uniformly bounded with $\Vert \frac{i e^{i\alpha}}{k} \mathscr{B}_+ \Vert \leq 6$,
then as in \eqref{eq5,6} it can be shown that
\begin{equation}\label{eq7,10}
\small{\Tilde{D}(k,r) = 
\exp \left( \mathcal{O} \left( \textup{Tr} \hspace{0.4mm} 
\one_{[\frac{r}{2},4r]} (\varepsilon\mathscr{B}_q) \right) \right).
}
\end{equation}
Now, establish a lower bound of $\Tilde{D}(ik,r)$ for $0 < \frac{2r}{3} < \vert k 
\vert < \frac{3r}{2}$, $k \in \br_+e^{-i\beta}$, $\beta > 0$ such that 
$z_q(ik) = 2bq - k^2 \in \Omega_{q,\nu}^{+} \Big( \frac{4r^2}{9},\frac{9r^2}{4} \Big)$,
$0 < \nu < \frac{8r^2}{9}$, is not a discrete eigenvalue of $H_\varepsilon$. Under 
this condition, thanks to \eqref{eq7,7} and \eqref{eq7,8},
$
I + \left( I + \frac{e^{i\alpha}}{k} \mathscr{B}_- + \varepsilon 
\mathscr{A}_{q}(k) \right)^{-1} \frac{e^{i\alpha}}{k} \mathscr{B}_+
$
is invertible. On the other hand, by exploiting the fact that
$\mathscr{B}_+\mathscr{B}_- = \mathscr{B}_-\mathscr{B}_+ = 0$, we get
\begin{equation}\label{eq7,11}
\begin{split}
\Big( I & + \frac{e^{i\alpha}}{k} \mathscr{B}_- + \varepsilon 
\mathscr{A}_{q}(k) \Big)^{-1} \frac{e^{i\alpha}}{k} \mathscr{B}_+ 
\\
& = \left[ I -
\left( I + \frac{e^{i\alpha}}{k} \mathscr{B}_- + \varepsilon 
\mathscr{A}_{q}(k) \right)^{-1} \left( \frac{e^{i\alpha}}{k} 
\mathscr{B}_- + \varepsilon \mathscr{A}_{q}(k) \right) \right]
\frac{e^{i\alpha}}{k} \mathscr{B}_+ \\
& = \frac{e^{i\alpha}}{k} \mathscr{B}_+ + \mathcal{O}(\varepsilon).
\end{split}
\end{equation}
Then, for $f \in L^{2}(\br^3) \setminus {\rm Ker}\, (\mathscr{B}_+)$,
we have
\begin{equation}\label{eq7,12}
\begin{split}
& \Biggl\vert \Im \left( \Big\langle \Big( I + \frac{e^{i\alpha}}{k} 
\mathscr{B}_- + \varepsilon \mathscr{A}_{q}(k) \Big)^{-1}
\frac{e^{i\alpha}}{k} \mathscr{B}_+ f,f \Big\rangle \right) \Biggr\vert \\
& = \Biggl\vert \Im \left( \Big\langle \Big( \frac{e^{i\alpha}}{k} 
\mathscr{B}_+ + \mathcal{O}(\varepsilon) \Big) f,f \Big\rangle 
\right) \Biggr\vert \\ 
& = \Bigl\vert \sin(\alpha + \beta) \Big\langle \frac{\mathscr{B}_+}{\vert k \vert} 
f,f \Big\rangle + \Im \left( \Big\langle \mathcal{O}(\varepsilon) 
f,f \Big\rangle \right) \Bigr\vert \geq \textup{Const.} \big\vert \sin(\alpha + \beta) 
\big\vert \Vert f \Vert^2,
\end{split}
\end{equation}
for $\varepsilon$ small enough and using the fact that $\sigma \big( \frac{1}{\vert k \vert} 
\mathscr{B}_+ \big) \subset \left] \frac{1}{3},6 \right[$. 
For $f \in {\rm Ker}\, (\mathscr{B}_+)$, we have
\begin{equation}\label{eq7,12}
\Re \left( \Big\langle \Big( I + \Big( I + \frac{e^{i\alpha}}{k} 
\mathscr{B}_- + \varepsilon \mathscr{A}_{q}(k) \Big)^{-1}
\frac{e^{i\alpha}}{k} \mathscr{B}_+ \Big) f,f \Big\rangle \right) 
=  \Vert f \Vert^2.
\end{equation}
Thus,
\begin{equation}\label{eq7,13}
\left\Vert \left( I + 
\left( I + \frac{e^{i\alpha}}{k} \mathscr{B}_- + \varepsilon 
\mathscr{A}_{q}(k) \right)^{-1} \frac{e^{i\alpha}}{k} 
\mathscr{B}_+ \right)^{-1} \right\Vert \leq C(\alpha,\beta),
\end{equation}
where $C(\alpha,\beta)$ is a constant depending on $\alpha$ and $\beta$. Consequently, 
as in \eqref{eq7,10}, it can be shown that
\begin{equation}\label{eq7,14}
\begin{split}
\Tilde{D}(ik,r)^{-1} = \textup{det} \Big\lbrace I -
& \Big( I + \frac{e^{i\alpha}}{k} \mathscr{B}_- + \varepsilon 
\mathscr{A}_{q}(k) \Big)^{-1} \frac{e^{i\alpha}}{k} \mathscr{B}_+ 
\\
& \Big[ I + \Big( I + \frac{e^{i\alpha}}{k} \mathscr{B}_- + 
\varepsilon \mathscr{A}_{q}(k) \Big)^{-1}\frac{e^{i\alpha}}{k} 
\mathscr{B}_+ \Big]^{-1} \Big\rbrace \\
& \leq
\exp \left( \mathcal{O} \left( \textup{Tr} \hspace{0.4mm} 
\one_{[\frac{r}{2},4r]} (\varepsilon\mathscr{B}_q) \right) \right).
\end{split}
\end{equation}
Namely,
\begin{equation}\label{eq7,15}
\small{\Tilde{D}(ik,r) \geq
\exp \left( -C \left( \textup{Tr} \hspace{0.4mm} 
\one_{[\frac{r}{2},4r]} (\varepsilon\mathscr{B}_q) \right) 
\right)},
\end{equation}
for some constant $C > 0$. We conclude as in the proof of Theorem \ref{theo2} by 
using the Jensen Lemma \ref{la,1}.

\medskip

\textbf{(iii)} Counted with their multiplicity, denote $(\mu_j)_j$ the decreasing 
sequence of the non-vanishing eigenvalues of the operator $P_p \textbf{\textup{W}} P_q$. 
Following \cite[Lemma 7]{bon}, there exits a constant $\nu > 0$ such that 
\begin{equation}\label{eq7,16}
\# \big\lbrace j : \mu_j - \mu_{j+1} > \nu \mu_j \big\rbrace 
= \infty.
\end{equation}
Since $\mathscr{B}_q$ and $P_p \textbf{\textup{W}} P_q$ have the same non-vanishing 
eigenvalues, then there exists a decreasing sequence of positive numbers 
$(r_\ell)_\ell$ with $r_{\ell} \searrow 0$, satisfying for any $\ell \in \bn$ 
(see Figure 7.1)
\begin{equation}\label{eq7,17}
\textup{dist} \big( r_\ell,\sigma(\mathscr{B}_q) \big) \geq 
\frac{\nu r_\ell}{2}.
\end{equation}
Moreover, for any $\ell \in \bn$, there exists a path $\Tilde{\Sigma}_{\ell} := 
\partial \omega_\ell$ (see Figure 7.1) with
\begin{equation}\label{eq10,18}
\omega_\ell := \big\lbrace \Tilde{k} \in \bc :
0 < \vert \Tilde{k} \vert < r_0 : \vert \Im(\Tilde{k}) \vert 
\leq \delta \Re(\Tilde{k}) : r_{\ell + 1} \leq \Re(\Tilde{k}) 
\leq r_\ell \big\rbrace,
\end{equation}
enclosing the eigenvalues of the operator $\mathscr{B}_q$ contained in 
$[r_{\ell + 1},r_\ell]$. 

\begin{figure}[h]\label{fig 4}
\begin{center}
\tikzstyle{+grisEncadre}=[dashed]
\tikzstyle{blancEncadre}=[fill=white!100]
\tikzstyle{grisEncadre}=[densely dotted]
\tikzstyle{dEncadre}=[dotted]

\begin{tikzpicture}[scale=1]

\draw [->] (0,-2) -- (0,4);
\draw [->] (-1.5,0) -- (6,0);

\draw (1.9,-0.95) -- (1.9,0.95) -- (3.5,1.75) -- (3.5,-1.75) -- cycle;

\draw [grisEncadre] (0,0) -- (1.9,0.95);
\draw [grisEncadre] (3.5,1.75) -- (5.5,2.75);

\node at (5.7,3) {$\Im(\Tilde{k}) = \delta \hspace*{0.08cm} \Re(\Tilde{k})$};
\node at (2.5,-1.5) {$\Tilde{\Sigma}_{\ell}$};

\draw [->] (1.5,-0.5) -- (1.86,-0.05);
\node at (1.4,-0.7) {$r_{\ell + 1}$};

\draw [->] (3.9,-0.5) -- (3.54,-0.04);
\node at (4,-0.7) {$r_{\ell}$};

\node at (3.3,0) {\tiny{$\bullet$}};
\node at (3.3,0.2) {\tiny{$\mu_{j}$}};

\node at (3.7,0) {\tiny{$\bullet$}};
\node at (3.9,0.2) {\tiny{$\mu_{j-1}$}};

\node at (1.6,0) {\tiny{$\bullet$}};
\node at (2.2,0) {\tiny{$\bullet$}};

\node at (2.6,0) {\tiny{$\bullet$}};
\node at (2.7,0.2) {\tiny{$\mu_{j+1}$}};

\node at (1.2,0) {\tiny{$\bullet$}};
\node at (1,0) {\tiny{$\bullet$}};
\node at (0.5,0) {\tiny{$\bullet$}};

\node at (4.3,0) {\tiny{$\bullet$}};
\node at (4.6,0) {\tiny{$\bullet$}};

\end{tikzpicture}
\caption{Representation of the path $\Tilde{\Sigma}_\ell  
= \partial \omega_\ell$.}
\end{center}
\end{figure}
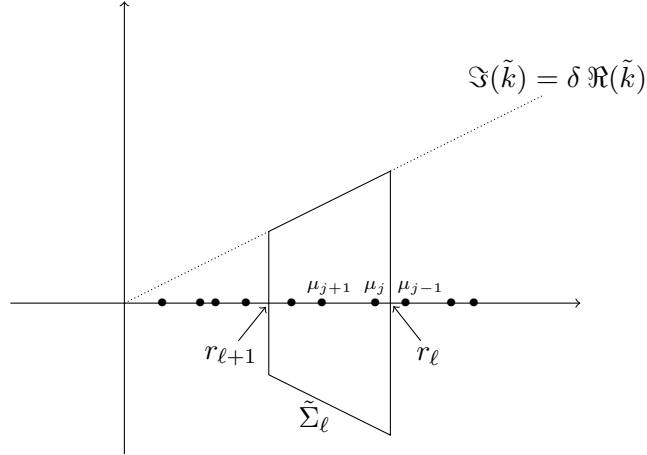

\noindent
It can be checked that the invertible operator $\Tilde{k} - \mathscr{B}_q$ 
for $\Tilde{k} \in \Tilde{\Sigma}_{\ell}$ satisfies 
\begin{equation}\label{eq7,19}
\big\Vert (\Tilde{k} - \mathscr{B}_q)^{-1} \big\Vert \leq 
\frac{\max \left( \delta^{-1}\sqrt{1 + \delta^2},(\nu/2)^{-1}\sqrt{1 + \delta^2} 
\right)}{ \vert \Tilde{k} \vert},
\end{equation}
uniformly with respect to $\Tilde{k} \in \Tilde{\Sigma}_{\ell}$. Introduce the path 
$\Sigma_{\ell} := -i\varepsilon e^{i\alpha} \Tilde{\Sigma}_{\ell}$ and estimate from 
below the number of zeros of $\textup{det}_{\lceil p \rceil} \big( I + 
\frac{i\varepsilon e^{i\alpha}}{k} \mathscr{B}_{q} + \varepsilon \mathscr{A}_{q}(k) 
\big)$ enclosed in $\big\lbrace z_q(k) \in \Omega_q^+(0,\eta^2) : k \in \omega_{\ell} 
\big\rbrace$, counted with their multiplicity. It is easy to see that according to 
the construction of $\Sigma_{\ell}$ and \eqref{eq7,19}, $I + \frac{i\varepsilon 
e^{i\alpha}}{k} \mathscr{B}_{q}$ is invertible for $k \in \Sigma_{\ell}$ and satisfies
\begin{equation}\label{eq7,20}
\left\Vert \left( I + \frac{i\varepsilon e^{i\alpha}}{k} 
\mathscr{B}_{q} \right)^{-1} \right\Vert \leq \max \left( 
\delta^{-1}\sqrt{1 + \delta^2}, (\nu/2)^{-1}\sqrt{1 + \delta^2} \right),
\end{equation}
uniformly with respect to $k \in \Sigma_{\ell}$. Then, for $k \in \Sigma_{\ell}$,
\begin{equation}\label{eq7,21}
\small{I + \frac{i \varepsilon e^{i\alpha}}{k} \mathscr{B}_q + \varepsilon \mathscr{A}_{q}(k) 
= \left( I + \varepsilon \mathscr{A}_{q}(k) \left( I + \frac{i \varepsilon e^{i\alpha}}{k} 
\mathscr{B}_q \right)^{-1} \right) \left( I + \frac{i \varepsilon e^{i\alpha}}{k} \mathscr{B}_q 
\right).}
\end{equation}
Choosing $0 < \varepsilon \leq \varepsilon_0$ small enough and using Property \textbf{g)} 
of Subsection \ref{ss3.1} given by \eqref{eq23,3}, we get for $k \in \Sigma_{\ell}$
\begin{equation}\label{eq7,22}
\small{\left\vert \textup{det}_{\lceil p \rceil} \left[ I + 
\varepsilon \mathscr{A}_{q}(k) 
\left( I + \frac{i \varepsilon e^{i\alpha}}{k} \mathscr{B}_q \right)^{-1} 
\right] - 1 \right\vert < 1
}.
\end{equation}
Consequently, by the Rouché Theorem, the number of zeros of 
$\textup{det}_{\lceil p \rceil} \big( I + \frac{i\varepsilon 
e^{i\alpha}}{k} \mathscr{B}_{q} + \varepsilon \mathscr{A}_{q}(k) 
\big)$ enclosed in $\big\lbrace z_q(k) \in \Omega_q^+(0,\eta^2) 
: k \in \omega_{\ell} \big\rbrace$ counted with their multiplicity,
is equal to that of $\textup{det}_{\lceil p \rceil} \big( I + 
\frac{i \varepsilon e^{i\alpha}}{k} \mathscr{B}_q \big)$ enclosed in 
$\big\lbrace z_q(k) \in \Omega_q^+(0,\eta^2) : k \in 
\omega_{\ell} \big\rbrace$ counted with their multiplicity. Thanks to \eqref{eq4,14}, 
this number is equal to $\textup{Tr} \hspace{0.4mm} \one_{[r_{\ell +1},r_\ell]} 
\big( P_{q} \textbf{\textup{W}} P_{q} \big)$. So, we get \eqref{lb1} since the 
zeros of $\textup{det}_{\lceil p \rceil} \big( I + \frac{i\varepsilon e^{i\alpha}}{k} \mathscr{B}_{q} + \varepsilon \mathscr{A}_{q}(k) \big)$ are the discrete eigenvalues of $\chh_\varepsilon$ with the same multiplicity, thanks to Proposition \ref{prop3,2} and 
Property \eqref{eqa,52} applied to \eqref{eq7,21}. The infiniteness of the number of 
the discrete eigenvalues claimed follows from the fact that the sequence $(r_\ell)_\ell$ 
is infinite tending to zero. The proof is complete.

\section{Proof of Theorem $\ref{theo6}$: Dominated accumulation}\label{s8}

The proof goes as that of item \textbf{(i)} of Theorem \ref{theo4}. 

Let the assumptions of Theorem \ref{theo6} hold. Then, for any $q \in \bn$, 
we have
\begin{equation}\label{eq7,01}
I + \mathcal{T}_{\varepsilon W} \big( z_{q}(k) \big) = I \pm \frac{i\varepsilon 
e^{i\alpha}}{k} \mathscr{B}_{q} + \varepsilon \mathscr{A}_{q}(k), \qquad
k \in \mathcal{D}_\pm^\ast(\eta).
\end{equation}
The operator $I \pm \frac{i\varepsilon e^{i\alpha}}{k} \mathscr{B}_{q}$ satisfies 
the bound \eqref{eq7,2} for $k \in e^{i\alpha} \mathcal{C}_\delta$, uniformly with 
respect to $0 < \vert k \vert < \eta$. Then,
\begin{equation}\label{eq7,03}
\small{I + \mathcal{T}_{\varepsilon W} \big( z_{q}(k) \big) = \big( I + A_\pm(k) \big) 
\left( I \pm \frac{i\varepsilon e^{i\alpha}}{k} \mathscr{B}_{q} \right)},
\end{equation}
with
\begin{equation}\label{eq7,04}
\small{A_\pm(k) := \varepsilon \mathscr{A}_{q} (k) \left( I \pm \frac{i\varepsilon 
e^{i\alpha}}{k} \mathscr{B}_{q} \right)^{-1}}.
\end{equation}
From Proposition \ref{prop4,1}, we deduce that there exists a constant $C > 0$ 
such that $\Vert \mathscr{A}_{q}(k) \Vert \leq C$ uniformly with respect to 
$0 \leq \vert k \vert \leq \eta$.
Then, for $0 < \varepsilon \leq \Tilde{\varepsilon}_0$ small enough, $I + 
\mathcal{T}_{\varepsilon W} \big( z_{q}(k) \big)$ is invertible for 
$k \in e^{i\alpha} \mathcal{C}_\delta$. Therefore, $z_q(k)$ is not a discrete 
eigenvalue, which proves the theorem.

%\newpage

\section{Appendix}\label{sa}

In this Appendix, we recall the notion of the index (with respect to a positively 
oriented contour) of a holomorphic function and a finite meromorphic operator-valued 
function, see for instance \cite[Definition 2.1]{bo}. 

For $f$ a holomorphic function in a neighbourhood of a contour $\gamma$, the index 
of $f$ with respect to $\gamma$ is defined by 
\begin{equation}\label{eqa,5}
ind_{\gamma} \hspace{0.5mm} f 
:= \frac{1}{2i\pi} \int_{\gamma} \frac{f'(z)}{f(z)} dz.
\end{equation}
Noting that if $f$ is holomorphic in a domain $\Omega$ with $\partial \Omega = 
\gamma$, then the residues theorem implies that $\textup{ind}_{\gamma} \hspace{0.5mm} f$ 
coincides with the number of zeros of $f$ in $\Omega$, counted with their multiplicity. 

Consider $D \subseteq \mathbb{C}$ a connected open set, $Z \subset D$ being a pure point and 
closed subset, and $A : \overline{D} \backslash Z \longrightarrow \textup{GL}(\mathscr{H})$ 
(the class of invertible operators on $\mathscr{H}$) being a finite meromorphic operator-valued 
function and Fredholm at each point of $Z$. The index of $A$ with respect to the contour 
$\partial \Omega$ is defined by 
\begin{equation}\label{eqa,51}
\small{Ind_{\partial \Omega} \hspace{0.5mm} A := \frac{1}{2i\pi} 
\textup{Tr} \int_{\partial \Omega} A'(z)A(z)^{-1} dz 
= \frac{1}{2i\pi} \textup{Tr} \int_{\partial \Omega} A(z)^{-1} 
A'(z) dz}.
\end{equation} 
We have the following properties: 
\begin{equation}\label{eqa,52}
Ind_{\partial \Omega} \hspace{0.5mm} A_{1} A_{2} = Ind_{\partial 
\Omega} \hspace{0.5mm} A_{1} + Ind_{\partial \Omega} \hspace{0.5mm} 
A_{2},
\end{equation} 
and if $K(z)$ lies in the trace class operator, then 
\begin{equation}\label{eqa,53}
Ind_{\partial \Omega} \hspace{0.5mm} (I+K)= ind_{\partial \Omega} 
\hspace{0.5mm} \det \hspace{0.5mm} (I + K).
\end{equation} 
For more details, see \cite[Chap. 4]{goh}.

The following lemma contains a version of the well-known Jensen 
inequality, see for instance \cite[Lemma 6]{bon} for a proof.

\begin{lem}\label{la,1} 
Let $\Delta$ be a simply connected sub-domain of $\mathbb{C}$ 
and let $g$ be a holomorphic function in $\Delta$ with continuous 
extension to $\overline{\Delta}$. Assume that there exists 
$\lambda_{0} \in \Delta$ such that $g(\lambda_{0}) \neq 0$ and 
$g(\lambda) \neq 0$ for $\lambda\in \partial \Delta$, the boundary 
of $\Delta$. Let $\lambda_{1}, \lambda_{2}, \ldots, \lambda_{N} 
\in \Delta$ be the zeros of $g$ repeated according to their 
multiplicity. For any domain $\Delta' \Subset \Delta$, there exists 
$C' > 0$ such that $N(\Delta',g)$, the number of zeros $\lambda_{j}$ 
of $g$ contained in $\Delta'$, satisfies
\begin{equation}\label{eqa,1}
N(\Delta',g) \leq C' \left( \int_{\partial \Delta} \textup{ln} 
\vert g(\lambda) \vert d\lambda - \textup{ln} \vert g(\lambda_{0}) 
\vert  \right).
\end{equation}
\end{lem}

%\textbf{Acknowledgements}.
%This work is partially supported by the ANR project NOSEVOL. The author is grateful to 
%J.-F Bony and V. Bruneau for suggest him this study. The author thanks V. Bruneau for 
%valuable discussion on the subject of this article.

%\newpage

\end{document}